\documentclass[fleqno,12pt,a4paper]{smfart}

\usepackage[all]{xy}
\usepackage[frenchb]{babel}
\usepackage{amsthm}
\usepackage{amsmath,amscd,amssymb,verbatim}
\usepackage{amsfonts}
\usepackage{amscd}

\renewcommand{\theequation}%
				{\thesubsection.\arabic{equation}}

\title{Rationalit\'e des s\'eries de Poincar\'e\\
et des\\
Fonctions Z\^eta motiviques}
\author{Julien Sebag}
\address{École normale supérieure, Département de mathématiques et applications, 
45 rue d'Ulm, 
75230 Paris cedex 05, France 
(UMR 8553 du CNRS)}
\email{Julien.Sebag@ens.fr}

\swapnumbers

\def\limproj{\mathop{\oalign{$\mathrm{lim}$\cr
\hidewidth$\longleftarrow$\hidewidth\cr}}}

\def\Spec{\mathop{\mathrm{Spec}}}
\def\Spf{\mathop{\mathrm{Spf}}}

\def\dim{{\mathrm{dim}}}
\def\max{\mathop{\mathrm{max}}}

\def\mod{\mathop{\mathrm{modulo}}}

\def\ord{{\mathrm{ord}}}

\def\G{{\mathrm{Gr}}}

\let\!=\mathcal

\newtheorem{defi}[subsubsection]{D\'efinition}
\newtheorem{lem}[subsubsection]{Lemme}
\newtheorem{cor}[subsubsection]{Corollaire}
\newtheorem{prop}[subsubsection]{Proposition}
\newtheorem{thm}[subsubsection]{Th\'eor\`eme}

\theoremstyle{remark}

\theoremstyle{remark}
\newtheorem{nota}[subsubsection]{Notations}

\theoremstyle{remark}
\newtheorem{exe}[subsubsection]{Exemple}

\theoremstyle{remark}
\newtheorem{exes}[subsubsection]{Exemples}

\theoremstyle{remark}
\newtheorem{rem}[subsubsection]{Remarque}

\theoremstyle{remark}

\begin{document}

\maketitle

\begin{abstract}
Soient $k$ un corps de caract\'eristique 0 et $X$ une $k[[t]]$-vari\'et\'e (\'eventuellement singuli\`ere) plate, purement de dimension relative $d$. Nous prouvons la rationalit\'e des s\'eries de Poincar\'e motiviques et de fonctions Z\^eta d'Igusa motiviques \`a l'aide de l'int\'egration motivique, du th\'eor\`eme de d\'esingularisation plong\'ee d'Hironaka, de la th\'eorie des mod\`eles de N\'eron faibles pour les sch\'emas formels et d'un th\'eor\`eme d'\'elimination des quantificateurs en th\'eorie des mod\`eles.
\end{abstract}

\section{Introduction}

\subsection{}
Dans \cite{dl1}, et par analogie avec le cas $p$-adique, Denef et Loeser ont d\'efini la notion de s\'eries de Poincar\'e motiviques et de fonctions Z\^eta motiviques pour une vari\'et\'e alg\'ebrique lisse sur un corps de caract\'eristique 0.
La g\'en\'eralisation de la th\'eorie de l'int\'egration motivique au cadre des sch\'emas formels sur un trait formel, \'etablie dans \cite{s}, permet d'\'etendre ces d\'efinitions au cas d'une vari\'et\'e $X$ (non suppos\'ee lisse) sur un anneau de s\'eries formelles \`a coefficients dans un corps $k$. Si $k$ est de caract\'eristique 0 (cette hypoth\`ese se justifie par l'utilisation du th\'eor\`eme de d\'esingularisation d'Hironaka et de th\'eor\`emes d'\'elimination des quantificateurs), ces s\'eries motiviques sont rationnelles, en un sens que l'on pr\'ecisera. 

\subsection{}
L'originalit\'e de cet article r\'eside essentiellement dans les th\'eor\`emes \ref{thm2} et \ref{thm3} qui apportent une r\'eponse explicite \`a la conjecture propos\'ee par Du Sautoy et Loeser (cf. remarque 3.11 de \cite{sl}). En ce qui concerne les fonctions Z\^eta motiviques, des r\'esultats analogues (cf. remarque \ref{looijenga}) \`a ceux que nous prouvons ont \'et\'e d\'emontr\'es dans \cite{lo} (cf. proposition 4.2 de \textit{loc. cit.}) par Looijenga, qui utilise une th\'eorie de l'int\'egration motivique pour les vari\'et\'es alg\'ebriques sur $k[[t]]$, avec $k$ de caract\'eristique 0. Il est important de noter que cette th\'eorie est compatible avec celle de \cite{s}, en remarquant que le sch\'ema des arcs d'une $k[[t]]$-vari\'et\'e alg\'ebrique est canoniquement isomorphe au sch\'ema de Greenberg de son compl\'et\'e $t$-adique. Le cadre des sch\'emas formels appara\^\i t donc comme le `` bon '' cadre pour parler d'int\'egrales motiviques.

\subsection{}
L'obstacle principal \`a la g\'en\'eralisation des diff\'erents r\'esultats au cadre des sch\'emas formels sur un trait formel est l'abscence (\`a notre connaissance) d'\'enonc\'es (g\'en\'eraux) sur les d\'esingularisations plong\'ees dans cette cat\'egorie. Remarquons que m\^eme si l'on consid\`ere des vari\'et\'es alg\'ebriques sur $k[[t]]$, avec $k$ de caract\'eristique 0, le th\'eor\`eme d'Hironaka fournit une r\'esolution `` r\'eguli\`ere '', \textit{i.e.} un $k[[t]]$-morphisme de sch\'emas $h~: Y\rightarrow X$, avec $Y$ r\'egulier, et non une r\'esolution `` lisse '', comme dans le cas trait\'e par Denef et Loeser. Pour palier cette carence, nous introduisons une classe de morphismes, que nous appelons \textit{r\'esolutions de N\'eron (faibles) plong\'ees}, sugg\'er\'ee par la th\'eorie des mod\`eles de N\'eron (cf. \cite{ls}). De tels morphismes $h~:Y\rightarrow X$ existent dans la cat\'egorie des $k[[t]]$-vari\'et\'es alg\'ebriques (cf. proposition \ref{existence desing Neron}) et r\'etablissent la condition de lissit\'e sur $Y$. En particulier, l'existence de tels morphismes permet d'affirmer qu'une $k[[t]]$-vari\'et\'e lisse v\'erifie les hypoth\`eses de la proposition 4.2 de \cite{lo} (cf. remarque \ref{looijenga}). Nous nous contentons donc d'exprimer les r\'esultats dans cette cat\'egorie, mais il appara\^\i tra \'evident au lecteur que les d\'emonstrations s'\'etendent au cas formel, toutes les fois que de tels morphismes existent.

\subsection{}
Le th\'eor\`eme de changement de variables pour l'int\'egrale motivique et l'existence de tels morphismes vont permettre de ramener l'\'etude de la rationalit\'e des diff\'erentes fonctions \`a l'\'evaluation dans le cas lisse de s\'eries associ\'ees \`a un diviseur \`a croisements normaux. Nous \'etudions donc au paragraphe 3 de telles s\'eries. Nous y montrons \'egalement comment d\'eduire de ces calculs la rationalit\'e de fonctions Z\^eta d'Igusa motiviques. Ces calculs interviendront \'egalement dans la derni\`ere partie consacr\'ee aux s\'eries de Poincar\'e motiviques. Les th\'eor\`emes g\'en\'eraux \ref{thm2} et \ref{thm3} n\'ecessitent l'introduction de notions et l'utilisation de r\'esultats de th\'eorie des mod\`eles, comme, par exemple, des théorèmes d'\'elimination des quantificateurs (cf. \cite{P}) que nous utiliserons essentiellement sous la forme du th\'eor\`eme 2.1 de \cite{dl1}. En outre, certains r\'esultats de base de cette th\'eorie appliqu\'ee en g\'eom\'etrie ont \'et\'e d\'emontr\'es dans \cite{dl1} et nous contenterons de renvoyer le lecteur aux preuves qui y sont faites. Par contre, nous allons montrer comment les ensembles d\'efinis \`a partir de \textit{conditions semi-alg\'ebriques} s'ins\`erent naturellement dans le langage de l'int\'egration motivique, et sont, comme on l'espérait, des ensembles mesurables au sens de \cite{s}. La rationalit\'e des s\'eries de Poincar\'e motiviques sera alors une simple cons\'equence du th\'eor\`eme \ref{thm2}.

\bigskip

Nous tenons \`a remercier Fran\c{c}ois Loeser pour les discussions que nous avons eues ensemble, qui ont vu na\^\i tre ce sujet sous cette forme g\'en\'erale et qui nous ont apport\'e une aide sensible. Nous souhaitons \'egalement remercier Antoine Chambert-Loir pour son aide et sa clairvoyance, qui, au travers de ses commentaires, nous ont permis d'am\'eliorer grandement une premi\`ere version de ce texte.

\section{Int\'egration motivique}

Soit $R$ un anneau de valuation discr\`ete complet, de corps r\'esiduel $k$ parfait et de corps des fractions $K$. On note $R_n:=R/(\pi)^{n+1}$, pour une uniformisante $\pi$ de $R$. Si $F$ est un corps parfait contenant $k$, on note $R_F:=R\hat{\otimes}_{L(k)} L(F)$, o\`u $L$ est le foncteur de la cat\'egorie des $k$-alg\`ebres dans celle des anneaux topologiques qui \`a $A$ associe~: $A$ dans le cas o\`u $R$ est d'\'egale caract\'eristique; $W(A)$ dans le cas o\`u $R$ est d'in\'egale caract\'eristique.

\subsection{Les sch\'emas formels}

\subsubsection{}

Dans ce texte, nous ne consid\`ererons que des sch\'emas formels \textit{$\pi$-adiques} et nous renvoyons \`a cite{s} ou \`a \cite{ega} \S10 pour de plus amples d\'etails. Un \textit{sch\'ema formel} $\!X$ est ainsi la donn\'ee d'un syst\`eme inductif de $R_n$-sch\'emas $(X_n)_{n\in\mathbf{N}}$ dont les morphismes de transition sont des immersions nilpotentes, induisant des isomorphismes $X_n\simeq X_{n+1}\otimes_{R_{n+1}}R_n$. En tant qu'espaces topologiques $\!X=X_0$ et l'on a un isomorphisme de faisceaux $\!O_{\!X}\simeq\limproj\!O_{X_n}$. Concr\`etement, la donn\'ee d'un sch\'ema formel affine $\!X=\Spf A$ est \'equivalente \`a celle d'une $R$-alg\`ebre $A$, munie de la topologie $\pi$-adique et compl\`ete pour celle-ci. Dans ce cas, $X_n$ est simplement le $R_n$-sch\'ema affine $\Spec (A\otimes_R R_n)$. On d\'esigne par $\mathbb{D}$ et $\mathbb{B}_R^N$ les sch\'emas formels $\Spf R$ et $\Spf R\{x_1,\ldots,x_N\}$ respectivement. Un \textit{morphisme de sch\'emas formels} $f:\!Y\rightarrow\!X$ est, dans notre contexte, la donn\'ee d'un syst\`eme inductif de morphismes de sch\'emas $(f_n)_{n\in\mathbf{N}}$, o\`u $f_n:Y_n\rightarrow X_n$. On notera $\underline{Form}_{/_\mathbb{D}}$ cette cat\'egorie. La sous-cat\'egorie des $\mathbb{D}$-sch\'emas formels \textit{sttf} est la sous-cat\'egorie pleine de $\underline{Form}_{/_\mathbb{D}}$ dont les objets $\!X$ sont s\'epar\'es (\textit{i.e.} tels que le sch\'ema $X_0$ soit s\'epar\'e) et topologiquement de type fini (\textit{i.e.} tels que $\!X$ admette un recouvrement \textit{fini} par des sch\'emas formels affines de la forme $\Spf A$, o\`u $A$ est une $R$-alg\`ebre topologique isomorphe \`a un quotient de la $R$-alg\`ebre topologique $R\{x_1,\ldots,x_N\}$).

\subsubsection{}

Dans cet article, nous consid\'erons essentiellement les compl\'et\'es $\pi$-adiques de $R$-vari\'et\'es alg\'ebriques, \textit{i.e.} les $\mathbb{D}$-sch\'emas formels \textit{sttf} alg\'ebrisables (cf. \cite{ega} \S10.8 pour la th\'eorie g\'en\'erale). Si $X$ est une $R$-sch\'ema (usuel), on note $\widehat{X}$ son compl\'et\'e $\pi$-adique. Remarquons que, si $X$ est plat, son compl\'et\'e $\widehat{X}$ l'est aussi. Conform\'ement \`a \cite{s}, on appelle dimension relative du $\mathbb{D}$-sch\'ema formel $\widehat{X}$ la dimension de $X_0$ (qui est aussi celle de $\widehat{X}_K$); et l'on dit que $\widehat{X}$ est \'equimensionnel si $\widehat{X}_K$ est \'equidimensionnel (ce qui implique que $X_0$ l'est aussi).

\subsubsection{}

Nous adoptons d\'esormais les conventions d'\'ecriture suivantes~:

\begin{enumerate}

\item[\textsl{(i)}] une vari\'et\'e (usuelle) sur $\Spec R$, \textit{i.e.} un $R$-sch\'ema de type fini, r\'eduit et s\'epar\'e (\textit{non suppos\'e irr\'eductible}), sera d\'esign\'ee par l'usage des caract\`eres romains \textit{italiques}~: $X$. La fibre sp\'eciale d'une $R$-vari\'et\'e $X$, \textit{i.e.} le $k$-sch\'ema $X\times_R k$, sera d\'esign\'ee par l'usage des caract\`eres \textit{lin\'eaux}~: $\mathsf{X}$.   

\item[\textsl{(i)}] un sch\'ema formel \textit{sttf} sur $\mathbb{D}$ sera d\'esign\'e par l'usage des caract\`eres \textit{calligraphiques}~: $\!X$. La fibre sp\'eciale d'un $\mathbb{D}$-sch\'ema formel, \textit{i.e.} le $k$-sch\'ema $X_0$, sera not\'ee $\mathsf{X}$.   

\end{enumerate}

\subsection{Le foncteur de Greenberg}
\label{Greenberg}

Nous renvoyons le lecteur \`a cite{s} ou \cite{ls} pour les d\'etails (cf. \'egalement les pages 276-277 de \cite{blr}). Pour tout $n\in\mathbf{N}$, on d\'efinit le foncteur de Greenberg de la cat\'egorie des $\mathbb{D}$-sch\'emas formels \textit{sttf} dans la cat\'egorie des $k$-sch\'emas s\'epar\'es, de type fini, qui \`a $\!X\in\underline{{Form}}_{/_\mathbb{D}}$ associe le $k$-sch\'ema r\'eduit $\G_n(\!X)$, dont les $F$-points, pour tout corps $F$ parfait contenant $k$, sont les $\big(R_F/(\pi)^{n+1}$\big)-points de $X_n$. Les $k$-vari\'et\'es $\G_n(\!X)$ forment un syst\`eme projectif admettant une limite (non localement de type fini en g\'en\'eral) dans la cat\'egorie des $k$-sch\'emas, que l'on note $\G(\!X)$ et qui est s\'epar\'ee et r\'eduite. Les $F$-points de $\G(\!X)$ sont les sections du $\big(\Spf R_F\big)$-sch\'ema formel $\!X_F:=\!X\hat{\times}_{\mathbb{D}} \Spf R_F$. 

Concr\`etement, si $\!X$ est affine d\'efini par les \'equations $f_i\in R\{X_1,\ldots,X_N\}$ pour $1\leq i\leq m$, la donn\'ee d'un point rationnel de $x\in\G(\!X)(k)$ est \'equivalente \`a celle d'un $N$-uplet $\varphi_x\in R^N$, v\'erifiant $f_i(\varphi_x)=0$ pour tout $1\le i\leq m$. En particulier, on en d\'eduit, pour tout $n\in\mathbf{N}$, des isomorphismes~:
$$
\G_n(\mathbb{B}_R^N)\simeq \mathbf{A}^{N(n+1)}_k.
$$

\begin{nota}
\label{remarque importante}
\begin{enumerate}

\item Pour tout $n\in\mathbf{N}$, les morphismes canoniques seront toujours notés de la manière suivante~:
$$
\xymatrix{\G(\!X)\ar[r]^{\pi_{n,\!X}}\ar[dr]_{\pi_{n-1,\!X}}& \G_n(\!X)\ar[d]^{\theta^n_{n-1}} \\
& \G_{n-1}(\!X).}
$$
Les morphismes $\pi_{n,\!X}$ (ou $\pi_n$) sont les morphismes de troncation. Les morphismes $\theta_n^{n-1}$ sont les morphismes de transition.

\item Soit $\!Y$ un $\mathbb{D}$-schéma formel \textit{sttf}. Soit $h:\!Y\rightarrow\!X$ un $\mathbb{D}$-morphisme de schémas formels. On notera encore $h$ le $k$-morphisme de schémas $\G(h)$, et $h_n$ le $k$-morphisme $\G_n(h)$, pour tout $n\in\mathbf{N}$. Ces morphismes $h$ et $h_n$ rendent commutatif le diagramme suivant~:
$$
\xymatrix{\G(\!Y)\ar[r]^{h}\ar[d]_{\pi_{n,\!Y}}&\G(\!X)\ar[d]^{\pi_{n,\!X}}\\
\G_n(\!Y)\ar[r]_{h_n}&\G_n(\!X)}
$$

\item Si $X$ est une $R$-vari\'et\'e, on note $\G(X)$ au lieu de $\G(\widehat{X})$ et $\pi_{n,X}$ (\textit{resp}. $\theta_{n,X}^{m}$) \`a la place de $\pi_{n,\widehat{X}}$ (\textit{resp}. $\theta_{n,\widehat{X}}^{m}$).

\end{enumerate}
\end{nota}

Soit $\!X$ un $\mathbb{D}$-sch\'ema formel \textit{sttf}, r\'eduit, plat, purement de dimension relative $d$. Soient $\!Z\hookrightarrow \!X$ un sous-$\mathbb{D}$-sch\'ema formel ferm\'e de $\!X$ et $\!I_\!Z\subset\!O_\!X$ son faisceau d'id\'eaux. Rappelons comment associer \`a tout point $x\in\G(\!X)$ sa multiplicit\'e le long de $\!Z$, qui est un entier positif, ou $+\infty$ si $x\in\G(\!Z)$. On d\'efinit l'application $\mathrm{mult}_{\!Z}~:\G(\!X)\rightarrow\mathbf{N}\cup\{\infty\}$ de la mani\`ere suivante~: soient $x\in\G(\!X)$ et $\varphi~:\Spf R_F\rightarrow\!X$ le $\mathbb{D}$-morphisme de sch\'emas formels correspondant \`a $\overline{x}$ par adjonction, avec $F\supset\kappa(x)\supset k$ un corps parfait et $\overline{x}\in\G(\!X)(F)$ de localit\'e $x$. Si le faisceau d'id\'eaux d\'efini comme l'image du morphisme canonique de $\!O_{R_F}$-modules $\varphi^\ast\!I_{\!Z}\rightarrow\!O_{R_F}$ est nul, on pose $\mathrm{mult}_{\!Z}(x)=\infty$; sinon, il existe un entier $n\in\mathbf{N}$ tel que $\varphi^\ast\!I_{\!Z}.\!O_{R_F}=\pi^n\!O_{R_F}$ et on pose $\mathrm{mult}_{\!Z}(x)=n$. Cet entier $n$ est \'egalement \'egal \`a la valuation du g\'en\'erateur de l'id\'eal $\!I_{\!Z,\pi_0(x)}.R_F$ et ne d\'epend pas du corps parfait $F$. Si $\!Z=\widehat{Z}$ et si $\!X=\hat{X}$, on note simplement $\mathrm{mult}_Z$ au lieu de $\mathrm{mult}_{\widehat{Z}}$. 

\begin{exes}
\begin{enumerate}

\item Soit $h:\!Y\rightarrow\!X$ un $\mathbb{D}$-morphisme de sch\'emas formels, tous deux purement de dimension relative $d$. Si $\!Z_h$ est le lieu singulier du morphisme $h$, \textit{i.e.} d\'efini par le faisceau d'id\'eaux $\mathrm{Fitt}_d(\Omega^1_{\!Y/\!X})$, l'application $\mathrm{mult}_{\!Z_h}$ est not\'ee $\ord_\pi(\mathrm{Jac})_h$ et appel\'ee \textit{ordre du jacobien de $h$}. Comme on va le voir, cette application joue un r\^ole crucial dans le th\'eor\`eme de changement de variables.

\item Si $\!Z$ est déterminé par un $\mathbb{D}$-morphisme de schémas formels $f:\!X\rightarrow\mathbb{B}_R^1$, \textit{i.e.} $\!Z=f^{-1}(0)$, on notera $\ord_\pi f$ l'application $\mathrm{mult}_{\!Z}$. 

\end{enumerate}
\end{exes}

Les assertions 1,2,3 du  lemme suivant sont d\'emontr\'ees dans \cite{s} (cf. lemme 5.1.2) et les points 4 et 5 d\'ecoulent directement de la d\'efinition~:

\begin{lem}
\label{mult1}
L'application $\mathrm{mult}_{\!Z}:\G(\!X)\rightarrow\mathbf{N}\cup\{\infty\}$ v\'erifie~:
\begin{enumerate}

\item si $\!U\hookrightarrow \!X$ est un ouvert de $\!X$, $\mathrm{mult}_{\!Z}(x)=\mathrm{mult}_{\!Z\cap\!U}(x)$ pour tout $x\in\G(\!U)$.

\item Pour la fibre en l'infini, on a l'\'egalit\'e
$$
(\mathrm{mult}_{\!Z})^{-1}(\infty)=\G(\!Z)=\G(\!Z_{\mathrm{red}})=(\mathrm{mult}_{\!Z_{\mathrm{red}}})^{-1}(\infty).
$$

\item $\pi_{0,\!X}\big(\mathrm{mult}^{-1}_{\!Z}(n)\big)\subset \mathsf{Z}$ pour tout $n\geq 1$. En outre, si $\!Z$ est lisse, l'inclusion précédente est une égalité.

\item Si $\!Z'\hookrightarrow\!X$ est un sous-$\mathbb{D}$-sch\'ema formel ferm\'e de $\!X$ contenu dans $\!Z$, alors, si $x\in\G(\!X)$, $\mathrm{mult}_{\!Z'}(x)\leq \mathrm{mult}_{\!Z}(x)$.

\item Supposons $\!Z\subset \mathsf{X}$. On a $\mathrm{mult}_{\!Z}^{-1}(0)=\G(\!X\backslash\!Z)$, $\mathrm{mult}_{\!Z}^{-1}(1)=\G(\!X)\backslash\G(\!X\backslash\!Z)$ et pour $n\geq 2$, $\mathrm{mult}_{\!Z}^{-1}(n)\not=\emptyset$ si et seulement si $n=1$ ou $n=0$. 

\end{enumerate}
\end{lem}

\begin{lem}
\label{mult2}
Soit $f:\!Y\rightarrow\!X$ un $\mathbb{D}$-morphisme de sch\'emas formels sttf. Soit $\!Z\hookrightarrow\!X$ un sous-$\mathbb{D}$-sch\'ema formel de $\!X$. Soit $n\in\mathbf{N}$. Alors~:
$$
f\big(\mathrm{mult}_{f^{-1}(\!Z)}^{-1}(n)\big)\subset\mathrm{mult}_{\!Z}^{-1}(n).
$$
En outre, si $\G(f)$ est bijectif, cette inclusion est une \'egalit\'e.
\end{lem}

\begin{proof}
Soit $x\in\mathrm{mult}_{f^{-1}(\!Z)}^{-1}(n)$. L'assertion 1 du lemme \ref{mult1} permet de supposer que $\!X$ et $\!Y$ sont des $\mathbb{D}$-sch\'emas formels \textit{ttf} affines. On note $\!X:=\Spf A$, $\!Y:=\Spf B$ et $\!I$ d\'esigne d\'esormais un id\'eal de $A$. Pour all\'eger les notations, et quitte à étendre les scalaires, on peut supposer que $x$ est un point rationnel de $\G(\!Y)$. Soit $\varphi:\mathbb{D}\rightarrow\!Y$ le $\mathbb{D}$-morphisme de sch\'emas formels lui correspondant par adjonction. Le sch\'ema formel $f^{-1}(\!Z)$ est d\'efini dans $\!Y$ par l'id\'eal de $B$, image de $f^{\ast}\!I$ dans $B$ par le morphisme canonique
$$
f^{\ast}\!I\rightarrow f^{\ast}\!O_{\!X}=B,
$$
induit par le morphisme d'inclusion $\!I\subset A$. La valeur de la fonction $\mathrm{mult}_{f^{-1}(\!Z)}$ en $x$ est alors d\'etermin\'ee par l'image de $\varphi^{\ast}f^{\ast}\!I$ dans $\varphi^{\ast}\!O_{\!Y}=B=\varphi^{\ast}f^{\ast}\!O_{\!X}$. Or cette image d\'etermine \'egalement la valeur de la fonction $\mathrm{mult}_{\!Z}^{-1}$ en $f(x)$.
\end{proof}

\subsection{Les anneaux de Grothendieck}

Soit $k$ un corps. On note $K_0(\underline{Var}_{/k})$ le groupe abélien engendré par les symboles $[S]$, pour $S$ une variété sur $k$ (\textit{i.e.} un $k$-schéma de type fini réduit et séparé), avec les relations $[S]=[S']$ si $S$ et $S'$ sont isomorphes et $[S]=[S'] +[S\backslash S']$ si $S'$ est une sous-variété fermée de $S$. Ce groupe possède une structure naturelle d'anneau, dont le produit est induit par le produit fibré. La classe de $\Spec k$ est l'élément neutre de cet anneau; on la note $\mathbf{1}$. En référence au \textit{motif} de Lefschetz, on note également $\mathbf{L}$ la classe de la droite affine dans $\!M$.

Soit $S$ une $k$-variété. Tout sous-ensemble constructible $C$ de $S$ peut s'écrire comme réunion finie disjointe de sous-$k$-variétés de $S$, \textit{i.e.} il existe des sous-$k$-variétés $(S_i)_{1\leq i\leq n}$ de $S$ telles que
$$
C=\sqcup_{1\leq i\leq n} S_i.
$$
À tout ensemble constructible $C$ de $S$, on peut donc associer naturellement une \textit{unique} classe $[C]$ dans $K_0(\underline{Var}_{/k})$ de sorte que, si $C$ et $C'$ sont deux ensembles constructibles d'une variété $S$, on a la relation $[C\cup C']=[C]+[C']-[C\cap C']$. 

On désigne par~:
$$
\!M:=K_0(\underline{Var}_{/k})[\mathbf{L}^{-1}]
$$ 
le localisé de $K_0(\underline{Var}_{/k})$ par rapport au système multiplicatif $\{\mathbf{1},\mathbf{L}, \mathbf{L}^2,\ldots\}$. Soient $F^m\!M$ le sous-groupe de $\!M$ engendré par les $[S]\mathbf{L}^{-i}$ tels que $\dim S-i\leq -m$, et $\widehat{\!M}$ le séparé complété de $\!M$ suivant cette filtration. On note $F^\bullet$ cette filtration et $\overline{\!M}$ l'image de $\!M$ dans $\widehat{\!M}$.

La filtration $F^\bullet$ définit une topologie métrisable sur $\widehat{\!M}$. L'application~:
$$
\xymatrix{\parallel\hskip 1mm\parallel:\widehat{\!M}\ar[r]&\mathbf{R}_{\geq 0}}
$$
définie par~:
$$
\lVert{a}\rVert=
\begin{cases}
2^{-n}& \text{si $a\in F^n\widehat{\!M}$ et $a\not\in F^{n+1}\widehat{\!M}$},\\
0 & \text{$a=0$}. 
\end{cases}
$$
est la norme induite par cette filtration. Elle munit $\widehat{\!M}$ d'une structure d'anneau normé non archimédien.

\subsection{La mesure motivique}

Soit $\!X$ un $\mathbb{D}$-sch\'ema formel \textit{sttf}, r\'eduit, plat et purement de  dimension relative $d$. Supposons en outre que la fibre g\'en\'erique de $\!X_{\mathrm{sing}}$ est de codimension au moins 1 dans $X_K$. Dans ce paragraphe on note $\!M$ au lieu de $\!M_k$. La mesure $\mu_{\!X}$, que nous allons utiliser, est d\'efinie (cf. \cite{s}) comme une application $\sigma$-additive sur les ensembles mesurables de $\G(\!X)$ et \`a valeurs dans $\!M$. Ces derniers sont construits comme des approximations de parties constructibles de $\G(\!X)$, appel\'ees \textit{cylindres}.

\subsubsection{}

On dit qu'une partie $C$ de $\G(\!X)$ est un ensemble cylindrique de rang $n$ de $\G(\!X)$ (ou plus simplement un $n$-cylindre de $\G(\!X)$), si $C=\pi_{n,\!X}^{-1}(E_n)$, o\`u $E_n$ d\'esigne une partie constructible de $\G_n(\!X)$. On dit que $C$ est un cylindre si $C$ est un cylindre  d'un certain rang $n$. On note $\mathbf{C}_\!X$ l'anneau bool\'een des cylindres de $\G(\!X)$. Dans cet anneau bool\'een, on peut distinguer l'id\'eal $\mathbf{C}_{0\!X}$ des cylindres $B$ stables, qui sont des cylindres de rang un certain entier naturel $n$ v\'erifiant que la restriction des morphismes de transition $\theta_m~:\pi_{m+1,\!X}(B)\rightarrow\pi_{m,\!X}(B)$ est une fibration localement triviale par morceaux de fibre $\mathbf{A}_k^{(m-n)d}$ pour tout $m\geq n$ (cf. d\'efinition 4.2.1). En particulier, ceci entra\^\i ne que~:
$$
[\pi_{m,\!X}(B)]\mathbf{L}^{-(m+1)d}=[\pi_{n,\!X}(B)]\mathbf{L}^{-(n+1)d}\in\!M
$$
pour tout $m\geq n$.

\begin{exe}
Soit $e\in\mathbf{N}$. Un exemple important de cylindre est le sous-ensemble $\G^{(e)}(\!X)$ de $\G(\!X)$ défini par~:
$$
\G^{(e)}(\!X):=\G(\!X)\backslash(\pi_{e,\!X}^{-1}(\G_e(\!X_{\mathrm{sing}}))).
$$
On en déduit que~:
$$
\G(\!X)=(\bigcup_{e\in\mathbf{N}}\G^{(e)}(\!X))\sqcup \G(\!X_{\mathrm{sing}}).
$$
En outre, le lemme 4.5.3 de \cite{s} assure que ce cylindre est stable.
\end{exe}

\subsubsection{}

L'ensemble des cylindres \'etant trop petit pour d\'ecrire les ph\'enom\`enes g\'en\'eraux, on est naturellement amen\'e \`a étendre cette notion par celle des parties \textit{mesurables}. On dit que $A$ est mesurable dans $\G(\!X)$ si, pour tout $\varepsilon>0$, il existe un ensemble $I$ au plus dénombrable et une famille de cylindres $(A_i)_{i\in I\cup\{0\}}$ tels que~:
\begin{enumerate}

\item[(i)]  $A\triangle A_0\subset\displaystyle\bigcup_{i\in I}A_i$, où $\triangle$ désigne la différence symétrique. 

\item[(ii)]  $\lVert\mu_\!X(A_i)\rVert< \varepsilon$, pour tout $i\in I$.
\end{enumerate}
Dans ce cas, on dit qu'une telle famille de cylindres $(A_i)_{i\in I\cup\{0\}}$, que l'on note $(A_0;(A_i)_{i\in I})$, est une $\varepsilon$-approximation cylindrique de $A$. Le cylindre $A_0$ est appelé la partie principale de l'approximation cylindrique. On dit que $A$ est fortement mesurable si de plus on peut choisir $A_0\subset A$, pour tout $\varepsilon>0$. Les parties mesurables forment un anneau bool\'een, qui contient celui des parties cylindriques. On le note $\mathbf{D}_\!X$.

\begin{thm}
Soit $\!X$ un $\mathbb{D}$-sch\'ema formel sttf, plat, r\'eduit et purement de dimension relative $d$. Il existe une unique mesure $\sigma$-additive sur $\mathbf{D}_\!X$ \`a valeurs dans $\widehat{\!M}$ telle que~:

\begin{enumerate}

\item pour tout $n$-cylindre stable $B_0\in\mathbf{C}_{0,\!X}$, 
$$
\mu_\!X(B_0)=[\pi_{n,\!X}(B_0)]\mathbf{L}^{-(n+1)d}\in\overline{\!M}.
$$

\item Pour $A$ et $B$ dans $\mathbf{D}_\!X$, $\lVert{\mu}_\!X(A\cup B)\rVert\leq \max(\lVert{\mu}_\!X(A)\rVert,\lVert{\mu}_\!X(B)\rVert)$. Si $A\subset B$, $\lVert{\mu}_\!X(A)\rVert\leq\lVert{\mu}_\!X(B)\rVert$.
 
\end{enumerate}

De plus, cette mesure v\'erifie les propri\'et\'es suivantes~:

\begin{enumerate}
\item[3.] Pour tout $n$-cylindre $B\in\mathbf{C}_\!X$,  
$$
\mu_\!X(B)=\lim_{e\rightarrow +\infty}\mu_{\!X}(B\cap\G^{(e)}(\!X))\in\widehat{\!M}.
$$

\item[4.] Pour tout ensemble mesurable $A\in\mathbf{D}_\!X$ et $\varepsilon>0$, soit $(A_0(\varepsilon),(A_i)_{i\in I})$ une $\varepsilon$-approximation cylindrique de $A$. Alors la limite~:
$$
\mu_\!X(A)={\displaystyle{\lim_{\varepsilon\rightarrow 0}}}~{\mu}_\!X(A_0(\varepsilon))\in\widehat{\!M},
$$
existe dans $\widehat{\!M}$. En outre, cette limite est ind\'ependante du choix des approximations cylindriques de $A$.
\end{enumerate}

\end{thm}

\subsection{L'int\'egrale motivique}
\label{int motivique}

Soit $A$ un ensemble mesurable de $\G(\!X)$. On dit qu'une application $\alpha~:A\rightarrow\mathbf{Z}\cup\{\infty\}$ est mesurable si, pour tout $n\in\mathbf{N}$, l'ensemble $\alpha^{-1}(n)\subset A$ est mesurable dans $\G(\!X)$. On dit que $\alpha$ est exponentiellement int\'egrable si cette application est mesurable et si la s\'erie 
$$
\sum_{n\in\mathbf{Z}}\mu_\!X(\alpha^{-1}(n))\mathbf{L}^{-n}
$$
converge dans $\widehat{\!M}$. On peut alors d\'efinir l'int\'egrale motivique d'une application $\alpha:A\rightarrow\mathbf{Z}\cup\{\infty\}$ exponentiellement int\'egrable par~:
$$
\boxed{\int_A\mathbf{L}^{-\alpha}d\mu:=\sum_{n\in\mathbf{Z}}\mu_\!X(\alpha^{-1}(n))\mathbf{L}^{-n}.}
$$
En guise d'exemple, rappelons cet \'enonc\'e (cf. \cite{s})~:
 
\begin{lem}
\label{mult3}
Soient $\!X$ un $\mathbb{D}$-sch\'ema formel sttf et $\!Z\hookrightarrow\!X$ un sous-$\mathbb{D}$-sch\'ema formel ferm\'e de $\!X$. Alors~:
\begin{enumerate}

\item pour tout $n\in\mathbf{N}$, la fibre $(\mathrm{mult}_{\!Z})^{-1}(n)$ est un $n$-cylindre de $\G(\!X)$.

\item Si $\G(\!Z)$ est un ensemble de mesure nulle dans $\G(\!X)$, l'application $\mathrm{mult}_{\!Z}$ est exponentiellement int\'egrable sur $\G(\!X)$.
\end{enumerate}
\end{lem}

Soit $h~:\!Y\rightarrow\!X$ un $\mathbb{D}$-morphisme de sch\'emas formels. On appelle lieu sauvage de $h$ l'ensemble 
$$
\Sigma_h= h^{-1}(\G(\!X_{\mathrm{sing}}))\cup\{y\in\G(\!Y)\mid \ord_\pi(\mathrm{Jac})_h(y)=\infty\},
$$
et on dit que $h$ est temp\'er\'e sur $B$ si l'intersection de $B$ et du lieu sauvage de $h$ est de mesure nulle dans $\G(\!Y)$.

\bigskip

Le th\'eor\`eme essentiel de cette th\'eorie de l'int\'egration est le th\'eor\`eme de changement de variables suivant.

\begin{thm}
\label{tcvg}
Soient $\!X$ et $\!Y$ deux $\mathbb{D}$-schémas formels sttf, réduits, plats et purement de dimension relative $d$. Supposons que $\!Y\rightarrow \mathbb{D}$ est lisse. Soient $A$ et $B$ des ensembles fortement mesurables de $\G(\!X)$ et $\G(\!Y)$ respectivement. Soit $h~:\!Y\rightarrow\!X$ un $\mathbb{D}$-morphisme, temp\'er\'e sur $B$, de sch\'emas formels qui induit une bijection entre $B$ et $A$. Alors, pour toute application exponentiellement int\'egrable $\alpha:A\rightarrow\mathbf{Z}\cup\{\infty\}$, l'application $B\rightarrow\mathbf{Z}\cup\{\infty\}$ qui \`a $y$ associe $\alpha(h(y))+\ord_\pi(\mathrm{Jac})(y)$ est exponentiellement int\'egrable et on a la formule:
$$
\int_A\mathbf{L}^{-\alpha}d\mu=\int_B\mathbf{L}^{-\alpha\circ h-\ord_\pi(\mathrm{Jac})_h}d\mu.
$$
\end{thm}

\subsection{Les r\'esolutions de N\'eron}
\label{desingularisation Neron}

Soit $\!X$ un $\mathbb{D}$-sch\'ema formel \textit{sttf}, purement de dimension relative $d$. On appelle \textit{r\'esolution de N\'eron de $\!X$} tout $\mathbb{D}$-morphisme de sch\'emas formels $h~:\!Y\rightarrow \!X$ v\'erifiant que~:

\begin{enumerate}

\item[\textsl{(i)}] $\!Y\rightarrow\mathbb{D}$ est un $\mathbb{D}$-sch\'ema formel \textit{sttf}, lisse, purement de dimension relative $d$;

\item[\textsl{(ii)}] le $\mathbb{D}$-morphisme de sch\'emas formels $h$ est temp\'er\'e;

\item[\textsl{(iii)}] le $k$-morphisme de sch\'emas induit $h~: \G(\!Y)\rightarrow \G(\!X)$ est bijectif.
\end{enumerate}

\begin{rem}
La th\'eorie des mod\`eles de N\'eron formelle (cf. \cite{bs} th\'eor\`eme 3.1) montre (cf. également \cite{ls} th\'eor\`eme 1.7.2) que si la fibre g\'en\'erique de $\!X$ est lisse en tant que $K$-espace analytique rigide, de tels morphismes existent. En outre, la th\'eorie d\'evelopp\'ee dans \textit{loc. cit.} d\'epend essentiellement de l'existence de ces r\'esolutions. Par ailleurs, la preuve du point \textsl{(v)} de la proposition \ref{existence desing Neron} ci-dessous montre que, si $X$ est une $k[[t]]$-vari\'et\'e alg\'ebrique r\'eguli\`ere, alors $X$ poss\`ede une r\'esolution de N\'eron, qui est donn\'ee par l'inclusion $j~:X_{\mathrm{lisse}}\rightarrow X$ du lieu lisse dans $X$.
\end{rem}

Soit $\!Z\hookrightarrow \!X$ un sous-$\mathbb{D}$-sch\'ema formel ferm\'e. On appelle \textit{r\'esolution de N\'eron plong\'ee de $(\!X,\!Z)$} tout $\mathbb{D}$-morphisme de sch\'emas formels $h~:\!Y\rightarrow \!X$ v\'erifiant que~:

\begin{enumerate}

\item[\textsl{(i)}] $\!Y\rightarrow\mathbb{D}$ est un $\mathbb{D}$-sch\'ema formel \textit{sttf}, lisse, purement de dimension relative $d$;

\item[\textsl{(ii)}] le $\mathbb{D}$-morphisme de sch\'emas formels $h$ est temp\'er\'e;

\item[\textsl{(iii)}] le $\mathbb{D}$-morphisme de sch\'emas formels induit $h~:\!Y\backslash h^{-1}(\!Z)\rightarrow \!X\backslash\!Z$ est injectif et \'etale;

\item[\textsl{(iv)}] le sous-$\mathbb{D}$-sch\'ema formel ferm\'e $h^{-1}\big(\G(\!Z)\big)$ est un ensemble de mesure nulle dans $\G(\!Y)$;

\item[\textsl{(v)}] le $k$-morphisme de sch\'emas induit $h~: \G(\!Y)\backslash h^{-1}\big(\G(\!Z)\big) \rightarrow \G(\!X)\backslash \G(\!Z)$ est bijectif;

\item[\textsl{(vi)}]  il existe un recouvrement de $\!Y$ par des sous-$\mathbb{D}$-sch\'emas formels ouverts et affines $\!V$ et, sur chaque $\!V$, un $\mathbb{D}$-morphisme de sch\'emas formels $\!V\rightarrow\mathbb{B}^d_R$ \'etale, induit par $d$ sections $z_1,\ldots, z_d$ au-dessus de $\!V$ telles que $h^{-1}(\!Z)$ est d\'efini, sur $\!V$, par l'annulation de l'une des $z_1^{n_1}\ldots z_d^{n_d}$, avec $n_i\geq 0$ pour tout $1\leq i\leq d$. 
\end{enumerate}

Si $\!X=\widehat{X}$ et $\!Z=\widehat{Z}$ sont des compl\'et\'es $\pi$-adiques de $R$-vari\'et\'es alg\'ebriques, on dira simplement que $(X,Z)$ admet une r\'esolution de N\'eron plong\'ee.

\begin{prop}
\label{existence desing Neron}
Soit $k$ un corps de caract\'eristique 0. Soient $X$ une $R$-vari\'et\'e alg\'ebrique plate, purement de dimension relative $d$ et $Z\hookrightarrow X$ un sous-$R$-sch\'ema ferm\'e de $X$ contenant l'ensemble des points singuliers de $X$, \textit{i.e.} l'ensemble des points o\`u le morphisme $X\rightarrow\Spec R$ n'est pas lisse. Si $R=k[[t]]$, alors $(X,Z)$ admet une r\'esolution de N\'eron plong\'ee.
\end{prop}

\begin{proof}
Soit $h'~:X'\rightarrow X$ une r\'esolution plong\'ee de $(X,Z)$, \textit{i.e.} un morphisme propre, birationnel, de source $X'$ r\'eguli\`ere et tel que $(h^{-1}(Z))_{\mathrm{red}}$ soit un diviseur \`a croisements normaux, qui existe par un th\'eor\`eme d'Hironaka (cf. \cite{g}). Notons $j~:U\hookrightarrow X'$ le lieu des points de $X'$ o\`u le morphisme de sch\'emas $X'\rightarrow\Spec R$ est lisse. Montrons alors que le morphisme induit par $h'\circ j~:U\rightarrow X$ par passage aux compl\'et\'es est une r\'esolution de N\'eron plong\'ee  de $(X,Z)$.

Les conditions (i) et (iii) d\'ecoulent directement de la d\'efinition et des propri\'et\'es du morphisme $h'$. La condition (ii) d\'ecoule de (iv), puisque $\Sigma_h\subset h^{-1}(Z)\cap U$. L'assertion (iv) se d\'eduit du fait que, par d\'efinition de $h'$, $(h')^{-1}(Z)$ est un diviseur de $X'$ et du lemme 4.4.2 de \cite{s}.

Pour prouver (v), il nous suffit de prouver que, pour tout corps parfait $F\supset k$, l'application canonique~:
$$
\Big(\G(X')\backslash (h')^{-1}\big(\G(Z)\big)\Big)(F)\rightarrow\Big(\G(X)\backslash \G(Z)\Big)(F)
$$
est une bijection. Remarquons en outre que, si $Y$ est une $k[[t]]$-vari\'et\'e, l'application naturelle 
$$
\widehat{Y}\big(\Spf(R_F)\big)\rightarrow Y(R_F)
$$
est bijective, pour toute extension $R_F=F[[t]]$. On est donc ramen\'e au cadre alg\'ebrique. Comme $h'$ est un isomorphisme au-dessus du complémentaire de $Z$, le critère valuatif de propreté assure que l'application induite par $h'$ ci-dessus est bijective. Il nous suffit alors de prouver que, pour de telles extensions $R_F\supset R$ , l'application naturelle~:
$$
{U}(R_F)\simeq {X'}(R_F)
$$ 
est bijective. Ceci revient \`a montrer que tout $R$-morphisme $\Spec R_F\rightarrow X'$ se factorise par le lieu lisse de $X'$. Si une telle propri\'et\'e est vraie pour le $R_F$-morphisme induit $\Spec R_F\rightarrow X'\times_R R_F$, elle l'est pour le $R$-morphisme $\Spec R_F\rightarrow X'$ (cf. \cite{ega44} proposition 17.7.1). Comme le morphisme $\Spec R_F\rightarrow \Spec R$ est r\'egulier (cf. \cite{ega42} d\'efinition 6.8.1), il d\'ecoule du diagramme cart\'esien de $R$-morphismes de sch\'emas suivant~:
$$
\xymatrix{X'\times_R R_F\ar[r]\ar[d]&R_F\ar[d]\\
				X'\ar[r]&R}
$$
et des propositions 6.8.2 et 6.5.2/(ii) de \cite{ega42} que le sch\'ema $X'\times_R R_F$ est encore r\'egulier. On peut donc supposer que $k=F$. Dans ce cas l'affirmation est d\'emontr\'ee dans la proposition 3.1/2 de \cite{blr}.

Pour (vi), soit $x\in U_0$. Il d\'ecoule du lemme \ref{pi para reg} qu'on peut trouver un syst\`eme r\'egulier de param\`etres de l'anneau local $\!O_{U,x}$ de la forme $(\pi,z_1,\ldots,z_{d})$. Par ailleurs, les images de ces sections locales dans le module des différentielles $\Omega^d_{\!O_{U_0,x}/_k}$ forment un système de générateurs. La propriété de changement de base du faisceau des diff\'erentielles relatives donne l'isomorphisme $\Omega^d_{\!O_{U_0,x}/_k}\simeq\Omega^d_{\!O_{U,x}/R}\otimes_R k$. Il découle donc du lemme de Nakayama que les images des $d$ sections $z_1,\ldots,z_{d}$ engendrent $\Omega^d_{\!O_{U,x}/R}$. En particulier, ces sections locales induisent un morphisme étale d'un ouvert affine $V$ dans $\mathbf{A}_R^d$. Les compl\'et\'es $t$-adiques des ouverts $V$ v\'erifient les propri\'et\'es souhait\'ees.
\end{proof}

\begin{lem}
\label{pi para reg}
Soient $R$ un anneau de valuation discr\`ete et $Y$ une $R$-vari\'et\'e lisse purement de dimension relative $d$. Soit $x\in Y_0$. Il existe un syst\`eme r\'egulier de param\`etres de $\!O_{Y,x}$ contenant $\pi$.
\end{lem}

\begin{proof}
Comme $Y$ est lisse, les anneaux locaux $\!O_{Y,x}$ et $\!O_{Y_0,x}\simeq \!O_{Y,x}\otimes_R k$ sont r\'eguliers. Soit $\overline{y_1},\ldots,\overline{y_d}$ un syst\`eme r\'egulier de param\`etres de l'anneau local $\!O_{Y_0,x}$. La suite $(\pi,y_1,\ldots,y_d)$ est alors un syst\`eme de param\`etres qui engendrent l'id\'eal maximal $x$ de $\!O_{Y,x}$, puisque~:
$$
\!O_{Y,x}/(\pi,y_1,\ldots,y_d)\simeq \!O_{Y_0,x}/(\overline{y_1},\ldots,\overline{y_d}).
$$ 
L'assertion d\'ecoule alors de la proposition 17.1.1 de \cite{ega41}.
\end{proof}

\begin{rem}
\label{rkmult2}
Soit $h~:Y\rightarrow X$ une résolution de Néron plongée de $(X,Z)$. Comme, pour tout $n\in\mathbf{N}$, $\mathrm{mult}_{Z}^{-1}(n)\subset \G(X)\backslash \G(Z)$, il découle du lemme \ref{mult2} que~: 
$$
h\big(\mathrm{mult}_{h^{-1}(Z)}^{-1}(n)\big)=\mathrm{mult}_{Z}^{-1}(n).
$$
Cette remarque est essentielle pour les calculs que nous effectuerons dans les différentes parties.
\end{rem}

\section{Familles et ensembles semi-alg\'ebriques}

Dans ce paragraphe, nous allons d\'efinir une classe particuli\`ere d'ensembles mesurables que l'on appelle \textit{ensembles semi-alg\'ebriques} et qui vont intervenir dans la d\'emonstration des diff\'erents r\'esultats.

\subsection{Les conditions semi-alg\'ebriques}

Fixons d\'esormais une cl\^oture alg\'ebrique $k^{\mathrm{alg}}\supset k$ de $k$.
Si $F$ est un corps alg\'ebriquement clos, et si $x\in F((t))$, on appelle \textit{composante angulaire de $x$}, que l'on note $\overline{ac}(x)$ le coefficient du terme de plus petit degr\'e, si $x\not=0$, et $0$ sinon. La valuation de $x$, not\'ee $\ord_t(x)$ est l'entier \'egal au plus petit degr\'e de $x$, si $x\not=0$, et $+\infty$ sinon. On \'etend alors la relation d'ordre usuelle sur $\mathbf{Z}$ \`a $\mathbf{Z}\cup\{+\infty\}$ de la mani\`ere habituelle et on adopte les conventions~:
\begin{enumerate}

\item[\textsl{(i)}] $(+\infty)+\ell=+\infty$, $\forall \ell\in\mathbf{Z}$;

\item[\textsl{(ii)}] $+\infty \equiv \ell~\mod~d$, $\forall \ell\in\mathbf{Z}$.
\end{enumerate}

Soient $x_1,\ldots, x_m$ des variables sur $k^{\mathrm{alg}}((t))$ et $\ell_1,\ldots,\ell_r$ des variables sur $\mathbf{Z}$. On appelle \textit{condition semi-alg\'ebrique} une combinaison bool\'eenne de conditions de la forme~:
\begin{flushleft}
$(SAL1)~~\ord_tf_1(x_1,\ldots,x_m)\geq \ord_tf_2(x_1,\ldots,x_m)+L(\ell_1,\ldots,\ell_r)$,\\
$(SAL2)~~\ord_tf_1(x_1,\ldots,x_m)\equiv L(\ell_1,\ldots,\ell_r) \ \ \textrm{mod $d$}$,\\
$(SAL3)~~g(\overline{ac}(f_1(x_1,\ldots,x_m)),\ldots, \overline{ac}(f_{s}(x_1,\ldots,x_m)))=0$,
\end{flushleft}
o\`u les $f_i$ sont des polyn\^omes \`a coefficients dans $k[[t]]$, $g$ est un polyn\^ome sur $k$ et $L\in\mathbf{Z}[X_1,\ldots,X_r]$. Les conditions semi-alg\'ebriques de la forme ci-dessus sont appel\'ees \textit{conditions semi-alg\'ebriques \'el\'ementaires}.

\begin{exe}
Si $f\in k[[t]][x_1,\ldots,x_m]$, la condition $f=0$ (et donc $f\not=0$) est une condition semi-alg\'ebrique.
\end{exe}

\subsection{Les familles semi-alg\'ebriques d'une vari\'et\'e}

Si $X$ est une $R$-vari\'et\'e et si $r\in\mathbf{N}$, on appelle \textit{famille semi-alg\'ebrique de $\G(X)$} une famille $(A_\ell)_{\ell\in\mathbf{N}^r}$ d'ensembles de $\G(X)$ tels qu'il existe un recouvrement de $X$ par un nombre fini d'ouverts affines $(U_i)_{i\in I}$ et pour chaque $U_i$, une condition semi-alg\'ebrique $\theta$ v\'erifiant que~:
$$
A_l\cap\G(U_i)=\{x\in\G(U_i)|\theta\big(h_1(\varphi_x),\ldots,h_m(\varphi_x);\ell\big)\},
$$ 
o\`u les $h_j$ sont des fonctions r\'eguli\`eres sur $U$ et o\`u, pour tout $x\in\G(X)$, $\varphi_x~:\Spf R_{\kappa(x)}\rightarrow \widehat{X}$ est le $\mathbb{D}$-morphisme de sch\'emas formels obtenu par adjonction et correspondant au point $x$. En particulier, $h_i(\varphi_x)$ s'identifie canoniquement \`a un \'el\'ement de $\kappa(x)[[t]]$. Une telle famille est alors \textit{d\'efinissable} sur chaque $U_i$ par les conditions semi-alg\'ebriques $\theta$ et les sections r\'eguli\`eres $h_i$, pour $1\leq i\leq m$.

\begin{rem}
\label{restriction}
Soit $(A_\ell)_{\ell\in\mathbf{N}^r}$ une famille semi-alg\'ebrique de $\G(X)$. Par d\'efinition, il existe un recouvrement fini de $X$ par des ouverts affines $(U_i)_{i\in I}$ v\'erifiant la propri\'et\'e de la d\'efinition. Pour chaque $i\in I$, on peut recouvrir $U_i$ par un nombre fini d'ouverts affines $(V_{i,j})_{j\in J}$. Le recouvrement $(V_{i,j})_{j\in J, i\in I}$ de $X$ v\'erifie encore la propri\'et\'e de l'\'enonc\'e. De m\^eme, il est clair que, si $U\hookrightarrow X$ est un ouvert de $X$ la famille $(A_\ell\cap\G(U))_{\ell\in\mathbf{N}^r}$ est semi-alg\'ebrique dans $\G(U)$.
\end{rem}

Un \textit{ensemble semi-alg\'ebrique} est une famille semi-alg\'ebrique avec $r=0$. Si $(A_\ell)_{\ell\in\mathbf{N}^r}$ est une famille semi-alg\'ebrique de $\G(X)$, alors, pour tout $\ell\in\mathbf{N}^r$, l'ensemble $A_\ell$ est semi-alg\'ebrique. On note $\mathbf{B}_{X}$ l'ensemble des parties semi-alg\'ebriques de $\G({X})$.

\begin{lem}
L'ensemble $\mathbf{B}_X$ est un anneau bool\'een, \textit{i.e.}
\begin{enumerate}

\item $\G(X)$ et $\emptyset$ appartiennent \`a cet ensemble;

\item $\mathbf{B}_X$ est stable par r\'eunion et intersection finies;

\item $\mathbf{B}_X$ est stable par passage au compl\'ementaire.
\end{enumerate}
\end{lem}

\begin{proof}
D\'ecoule directement du fait que les conditions semi-alg\'ebriques sont stables par de telles op\'erations.
\end{proof}

\begin{lem}
\label{image inverse}
Soient $X$ et $Y$ deux $R$-vari\'et\'es et $g~:Y\rightarrow X$ un $R$-morphisme de sch\'emas. Si $(A_\ell)_{\ell\in\mathbf{N}^r}$ est une famille semi-alg\'ebrique de $\G({X})$, la famille 
$$
(g^{-1}(A_\ell))_{\ell\in\mathbf{N}^r}
$$
est semi-alg\'ebrique dans $\G({Y})$. En particulier, si $A$ est un ensemble semi-alg\'ebrique de $\G({X})$, l'ensemble $g^{-1}(A)$ est un ensemble semi-alg\'ebrique de $\G({Y})$. 
\end{lem}

\begin{proof}
La seconde assertion d\'ecoule de la premi\`ere par d\'efinition. Il nous suffit donc de prouver la premi\`ere. Soit $A_\ell$, $\ell\in\mathbf{N}^r$, une famille semi-alg\'ebrique de $\G(X)$. Par d\'efinition, il existe donc un recouvrement (fini) de $X$ par des ouverts affines $(U_i)_{i\in I}$ v\'erifiant que~:
$$
A_\ell\cap\G({U_i})=\{x\in\G({U_i})|\theta\big(h_1^{(i)}(\varphi_x),\ldots,h_{m_i}^{(i)}(\varphi_x);l\big)\},
$$ 
Pour tout $i\in I$, soit $(V_{i,j})_{j\in J_i}$ un recouvrement (fini) de $g^{-1}(U_i)$ par des ouverts affines de $Y$. En particulier, les $V_{i,j}$, pour tout $i\in I$ et tout $j\in J$, recouvrent $Y$ et v\'erifient 
$$
g^{-1}(A_\ell)\cap \G(V_{i,j})=\{y\in\G(V_{i,j})|\theta\big((h_1^{(i)}\circ g)(\varphi_y),\ldots,(h_{m_i}^{(i)}\circ g)(\varphi_y);\ell\big)\},
$$
pour tout $\ell\in\mathbf{N}^r$.
\end{proof}

\subsubsection{}
\label{atlas}
Soit $(A_\ell)_{\ell\in\mathbf{N}^r}$ une famille semi-alg\'ebrique de $\G(X)$. Quitte \`a restreindre les $U_i$, on peut supposer (cf. remarque \ref{restriction}) que $U_i$ est immerg\'e dans un espace affine de dimension $n_i$, \textit{i.e.} est muni d'une immersion ferm\'ee $U_i\hookrightarrow\mathbf{A}^{m_i}_R$. Par d\'efinition, il existe une condition semi-alg\'ebrique $\theta$ et $m_i$ sections r\'eguli\`eres $h^{(i)}_j$ de $U_i$ telles que~:
$$
A_\ell\cap\G({U_i})=\{x\in\G({U_i})|\theta\big(h_1^{(i)}(\varphi_x),\ldots,h_{m_i}^{(i)}(\varphi_x);l\big)\}.
$$ 
Soit $(x^{(i)}_1,\ldots,x^{(i)}_{n_i})$ un syst\`eme de coordonn\'ees de $U_i$, \textit{i.e.} $n_i$ sections r\'eguli\`eres sur $U_i$ induites par le morphisme d'immersion. Pour tout $1\leq j\leq m_i$, les $h_j^{(i)}$ sont des polyn\^omes en les $x_s^{(i)}$, \textit{i.e.} il existe $m_i$ polyn\^omes $H^{(i)}_j\in k[[t]][Y_1,\ldots,Y_{n_i}]$ tels que~:
$$
h^{(i)}_j=H^{(i)}_j(x^{(i)}_1,\ldots,x^{(i)}_{n_i}),
$$
$1\leq j\leq m_i$. Soit $\theta':=\theta\circ H^{(i)}$ la condition semi-alg\'ebrique d\'efinie par $\theta(H^{(i)}_1,\ldots,H^{(i)}_{m_i})$ v\'erifie que~:
$$
A_\ell\cap\G({U_i})=\{x\in\G({U_i})|\theta'\big(x_1^{(i)}(\varphi_x),\ldots,x_{m_i}^{(i)}(\varphi_x);l\big)\}.
$$ 

\begin{defi}
Si $(A_\ell)_{\ell\in\mathbf{N}^r}$ est une famille semi-alg\'ebrique de $\G(X)$, on appelle carte semi-alg\'ebrique de la famille $(A_\ell)_{\ell\in\mathbf{N}^r}$ un ouvert affine $U$ de $X$, muni d'une immersion ferm\'ee dans un espace affine $U\hookrightarrow\mathbf{A}^{m}_R$, tels que~:
$$
A_\ell\cap\G({U})=\{x\in\G({U})|\theta\big(x_1(\varphi_x),\ldots,x_m(\varphi_x);l\big)\},
$$ 
o\`u les $x_i$ forment un syst\`eme de coordonn\'ees de $U$ et $\theta$ une condition semi-alg\'ebrique. La donn\'ee d'un recouvrement fini de $X$ par des cartes semi-alg\'ebriques de la famille $(A_\ell)_{\ell\in\mathbf{N}^r}$ est un atlas semi-alg\'ebrique de la famille $(A_\ell)_{\ell\in\mathbf{N}^r}$.  
\end{defi}

Il d\'ecoule de ce qui pr\'ec\`ede que de tels atlas existent pour toute famille semi-alg\'ebrique de $\G(X)$. En outre, une famille $(A_\ell)_{\ell\in\mathbf{N}^r}$ de $\G(X)$ est semi-alg\'ebrique si et seulement si elle poss\`ede un atlas semi-alg\'ebrique.

\begin{lem}
Soit $(A_\ell)_{\ell\in\mathbf{N}^r}$ une famille semi-alg\'ebrique de $\G(X)$. Alors la donn\'ee d'une carte semi-alg\'ebrique (resp. d'un atlas semi-alg\'ebrique) de $(A_\ell)_{\ell\in\mathbf{N}^r}$ ne d\'epend pas des plongements. 
\end{lem}

\begin{proof}
On est ramen\'e au cas d'un ensemble semi-alg\'ebrique $A$ de $\G(X)$. Soit $j_1~:U\hookrightarrow\mathbf{A}_R^m$ une carte semi-alg\'ebrique de $A$. Soit $j_2~:U\hookrightarrow\mathbf{A}_R^n$ un second plongement de $U$ dans un espace affine. Soit $h~: \mathbf{A}_R^n\rightarrow \mathbf{A}_R^m$ le $R$-morphisme de changement de coordonn\'ees. Avec $j_1$, il existe, par d\'efinition, une condition semi-alg\'ebrique $\theta$ d\'efinissant $A\cap\G(U)$. La condition semi-alg\'ebrique $\theta\circ h$, d\'efinie comme dans la preuve de la proposition \ref{image inverse}, fait de $U$ (muni de $j_2$) une carte semi-alg\'ebrique de $A$.
\end{proof}

\begin{prop}
\label{image directe}
Soient $X$ et $Y$ deux $R$-vari\'et\'es et $g~:Y\rightarrow X$ un $R$-morphisme de sch\'emas. Si $(A_\ell)_{\ell\in\mathbf{N}^r}$ est une famille semi-alg\'ebrique de $\G({Y})$, la famille $(g(A_\ell))_{\ell\in\mathbf{N}^r}$ est semi-alg\'ebrique dans $\G({X})$. En particulier, si $A$ est un ensemble semi-alg\'ebrique de $\G({Y})$, l'ensemble $g(A)$ est un ensemble semi-alg\'ebrique de $\G({X})$. 
\end{prop}

\begin{proof}
Soit $(A_\ell)_{\ell\in\mathbf{N}^r}$ une famille semi-alg\'ebrique de $\G(X)$. La propri\'et\'e \'etant locale en $X$, on peut supposer que $X$ est affine. En outre, comme les conditions semi-alg\'ebriques sont stables par r\'eunion finie, on peut supposer que $Y$ est carte semi-alg\'ebrique de la famille $(A_\ell)_{\ell\in\mathbf{N}^r}$. Par d\'efinition,
$$
g(A_\ell)=\left\{y\in\G(X)~\vert~\exists x\in\G(Y)~\textrm{tel que}~g(x)=y~\textrm{et}~x\in A_\ell\right\}.
$$
Comme $A_\ell$ est semi-alg\'ebrique, il d\'ecoule d'un th\'eor\`eme d'\'elimination des quantificateurs, d\^u \`a Pas (cf. \cite{P} ou \cite{dl1} th\'eor\`eme 2.1), que $g(A_\ell)$ est semi-alg\'ebrique dans $\G(X)$ et que la famille $(g(A_\ell))_{\ell\in\mathbf{N}^r}$ est semi-alg\'ebrique. 
\end{proof}

\subsubsection{}
\label{simple}

Soit $A$ un ensemble semi-alg\'ebrique de $\G(\widehat{X})$. On dit qu'une fonction
$$
\alpha:A\times\mathbf{Z}^n\rightarrow\mathbf{Z}\cup\{\infty\}
$$
est \textit{simple} si la famille de parties
$$
\{x\in\G(\widehat{X})~|~\alpha(x,l_1,\ldots,l_n)=l_{n+1}\},
$$
avec $(l_1,\ldots,l_{n+1})\in\mathbf{N}^{n+1}$, est une famille semi-alg\'ebrique de $\G(X)$.

\subsection{Les ensembles mesurables et les ensembles semi-alg\'ebriques}

\begin{prop}
Soit $X$ une $R$-vari\'et\'e. Soit $A\subset\G(X)$ un sous-ensemble semi-alg\'ebrique. Alors l'ensemble $\pi_{n,X}(A)\subset \G_n(X)$ est constructible.
\end{prop}

\begin{proof}
Soient $n\in\mathbf{N}$ et $A$ un ensemble semi-alg\'ebrique de $\G(X)$. Par d\'efinition de la constructibilit\'e, on peut supposer que $X$ est affine et est une carte semi-alg\'ebrique de $A$. Il existe donc un entier $m\in\mathbf{N}$ et une immersion ferm\'ee $\G(X)\hookrightarrow\G(\mathbf{A}^m_R)$. Comme
$$
\pi_{n,X}(A)=\pi_{n,\mathbf{A}_R^m}(A)\cap\G(X),
$$
on peut supposer en outre que $X=\mathbf{A}_R^m$. L'exemple \ref{ordjac} assure que l'ensemble $B:=\pi_{n,X}^{-1}\left(\pi_{n,X}(A)\right)$ est semi-alg\'ebrique puisque $A$ l'est.

Si $F\supset k$ est une extension de corps, et si $x\in\left(F[t]/(t^{n+1})\right)^m$, avec pour $1\leq i\leq m$,
$$
x_i=\sum_{j=1}^n a_j^{(i)}t^j,
$$
o\`u $a_j^{(i)}\in F$; on note $s(x)$ le $m$-uplet de $F[[t]]^m$ d\'efinie par les coordonn\'ees~:
$$
\left(s(x)_i\right):=\sum_{j=1}^n a_j^{(i)}t^j.
$$
Il est clair que $x\in\pi_{n,X}(A)$ si et seulement si $s(x)\in B$. On est donc ramen\'e \`a \'etudier l'ensemble des \'el\'ements de $B$ de la forme $s(x)$ ci-dessus. Or l'ensemble $B$ s'exprime comme une combinaison bool\'eenne (\textit{a priori} infinie) de cylindres (cf. exemple \ref{restriction}). Le r\'esultat d\'ecoule alors du fait que si $y\in B$ de la forme $s(x)$, ces r\'eunions sont \`a support fini et d'un th\'eor\`eme de Schappacher \cite{sch}, qui implique que, pour tout sous-$R$-sch\'ema ferm\'e $Z$ de $X$~:
$$
\pi_{n,X}\left(\G(Z)\right)=\theta^N_n\left(\G_N(Z))\right),
$$
pour un entier $N$ suffisamment grand.
\end{proof}

\begin{exes}
\label{ordjac}
\begin{itemize}
\item[1.] Si $A$ est un ensemble semi-alg\'ebrique de $\G(X)$, alors l'ensemble~:
$$
\pi_{n,X}^{-1}\left(\pi_{n,X}(A)\right)
$$
est semi-alg\'ebrique dans $\G(X)$. En effet, supposons que $X$ est une carte semi-alg\'ebrique de $A$. Il d\'ecoule de l'\'egalit\'e~:
$$
B:=\pi_{n,X}^{-1}\left(\pi_{n,X}(A)\right)=\left\{y\in\G(X)~\mid~\exists x\in A, y\equiv x\mod[t^{n+1}]\right\},
$$
et du th\'eor\`eme d'\'elimination des quantificateurs 2.1 de \cite{dl1}, que $B$ est un ensemble semi-alg\'ebrique de $\G(X)$.

\item[2.] Si $Z\hookrightarrow X$ est un sous-$R$-sch\'ema ferm\'e de $X$, alors la famille~:
$$
\left(\mathrm{mult}_Z^{-1}(n)\right)_{n\in\mathbf{N}}
$$
est une famille semi-alg\'ebrique de $\G(X)$. En effet, supposons que $X$ est une carte semi-alg\'ebrique pour cette famille. Soit $(z_i)_{1\leq i\leq s}$ une pr\'esentation de l'id\'eal $I\subset  R[y_1,\ldots y_N]$ d\'efinissant $Z$ dans $X$. La famille $\left(\pi_{n,X}^{-1}(\G_n(X))\right)_{n\in\mathbf{N}}$ est alors d\'efinissable par la condition semi-alg\'ebrique~:
$$
\theta(x_1,\ldots,x_s;n):=\left\{\begin{array}{l}
f_i(x_1,\ldots,x_s)=x_i,~~\forall 1\leq i\leq s\\
L(y)=y\\				
\ord_t f_i\geq L(n),~~\forall 1\leq i\leq s
\end{array}\right.,
$$
et les sections r\'eguli\`eres $z_i$ de $X$, pour $1\leq i\leq s$. L'assertion d\'ecoule alors de la relation~:
$$
\mathrm{mult}_Z^{-1}(n)=\pi_{n-1,X}^{-1}(\G_{n-1}(Z))\backslash\pi_{n,X}^{-1}(\G_n(Z))
$$
et des propri\'et\'es de stabilit\'e par intersection et passage au compl\'ementaire.
\end{itemize}
\end{exes}

\begin{exe}
La fonction $x\mapsto\mathrm{mult}_{Z}(x)$, o\`u $Z\hookrightarrow X$ est un sous-$R$-sch\'ema ferm\'e de $X$, est simple~: ceci est une simple retraduction de l'exemple \ref{ordjac}/2 ci-dessus.
\end{exe}

\subsubsection{}

Soit $A$ un ensemble semi-alg\'ebrique de $\G(X)$. On dit que $A$ est un \textit{ensemble faiblement stable de rang $n\in\mathbf{N}$} si $A$ est la r\'eunion de fibres du $k$-morphisme $\pi_{n,X}:\G(X)\rightarrow\G_n(X)$. On dit que $A$ est \textit{faiblement stable} s'il existe $n\in\mathbf{N}$ tel que $A$ soit faiblement stable de rang $n$.

\begin{lem}
\label{cylindre et ens stable} 
Les ensembles faiblement stables forment un anneau bool\'een, que l'on note $\mathbf{B}_{X}^{fs}$. En outre, si $A$ est un ensemble faiblement stable de rang $n$, alors $A$ est un $n$-cylindre.
\end{lem}

\subsubsection{}

Supposons que $X$ est purement de dimension relative $d$. Soit $A$ un ensemble semi-alg\'ebrique. On dit que $A$ est \textit{stable de rang $n\in\mathbf{N}$} si $A$ est faiblement stable au rang $n$ et si l'application \'evidente $\pi_{m+1,\widehat{X}}(\G(X))\rightarrow\pi_{m,\widehat{X}}(\G(X))$ est une fibration localement triviale au-dessus de $\pi_{m,\widehat{X}}(A)$ de fibre $\mathbf{A}_k^d$ pour tout $m\geq n$ (cf. d\'efinition 4.2.1 de \cite{s}). On dit que $A$ est \textit{stable} s'il existe $n\in\mathbf{N}$ tel que $A$ est stable de rang $n$. Il d\'ecoule des d\'efinitions et du lemme \ref{cylindre et ens stable} que tout ensemble stable est un cylindre stable. On note $\mathbf{B}_{X}^{s}$ l'ensemble des parties semi-alg\'ebriques stables.

\begin{lem}
\label{sal et mesurable}
Si $A\in\mathbf{B}_{X}$ est un ensemble semi-alg\'ebrique de $\G(X)$, alors il existe un ensemble de mesure nulle $M\subset\G(X)$ de $\G(X)$ et une famille semi-alg\'ebrique $(A_i)_{i\in\mathbf{N}}$ de $\G(X)$, avec $A_i\subset A$, tels que~:
\begin{enumerate}

\item les $A_i$ et $M\cap A$ forment une partition de $A$;

\item les $A_i$ sont des ensembles semi-alg\'ebriques stables de rang $n_i$;

\item  $\lim_{i\rightarrow+\infty}(\dim (\pi_{n_i,\widehat{X}}(A_i))-(n_i+1)d)=-\infty$.
\end{enumerate}
\end{lem}

\begin{proof}
On peut supposer que $X$ est une carte semi-alg\'ebrique de $A$ et est irr\'eductible. Soit $f_i$, pour $1\leq i\leq s$, les polyn\^omes non nuls \`a partir desquels $A$ est d\'efini. Soit $g$ une section r\'eguli\`ere de $X$ s'annulant sur le lieu singulier de $X$, \textit{i.e.} $X_{\mathrm{sing}}\subset V(g)\subsetneq X$. En particulier que $\widehat{X}_{\mathrm{sing}}\subset V(g)\subsetneq \widehat{X}$. Soit $F:=g\times\prod_i f_i$. Notons $S:=V(F)\subsetneq X$. Pour $i\in\mathbf{N}$, soit
$$
A_i:=\{x\in A~|~\ord_tF(x)=i\}.
$$
Par d\'efinition, $A_i\subset A$ est un ensemble semi-alg\'ebrique comme intersection de deux ensembles semi-alg\'ebriques (cf. exemple \ref{ordjac}/2). Soient $y\in\pi_{i,X}^{-1}(A_i)$ et $y\in k'[[t]]^N$, pour $k'\supset k$, la s\'erie formelle obtenue par adjonction. Il d\'ecoule d'un d\'eveloppement de Taylor que, pour tout $u\in k'[[t]]$,
$$
Q(y+t^{i+1}u)=Q(y)+t^{i+1}(\ldots),
$$
pour tout polyn\^ome $Q$ \`a coefficients dans $R$. Cette remarque entra\^ine en particulier que $A_i$ est faiblement stable. En outre, $A_i\subset\G(X)\backslash\G(S)\subset\G(X)\backslash\G(X_{\mathrm{sing}})$. Par suite, $A_i$ est contenu, par quasi-compacit\'e de la topologie constructible (cf. lemme 4.3.7 de \cite{s}), dans un nombre fini d'ensembles $\G^{(n(i))}(X)$. Comme la suite de ces ensembles est croissante, il existe $n_0(i)\in\mathbf{N}$ tel que $A_i\in\G^{(n_0(i))}(X)$.  Par cons\'equent, $A_i$ est stable. Le point 3 d\'ecoule alors du fait que la suite 
$$
(\dim(\pi_{i,\widehat{X}}((\ord_t F)^{-1}(i)))-(i+1)d)_{i\in\mathbf{N}}
$$
tend vers $-\infty$ par le lemme 4.4.2 de \cite{s}, puisque $S$ est de codimension 1. 
\end{proof}

\begin{prop}
\label{sal et mesurables}
Tout ensemble semi-alg\'ebrique de $\G(X)$ est mesurable, \textit{i.e.} $\mathbf{B}_{X}\subset \mathbf{D}_{\widehat{X}}$. En particulier, si $(A_\ell)_{\ell\in\mathbf{N}^r}$ est une famille semi-alg\'ebrique de $\G(X)$ et si $\alpha~:A\cup\mathbf{N}^r\rightarrow \mathbf{Z}\cup\{\infty\}$ est une fonction simple, alors, pour tout $\ell\in\mathbf{N}^r$, l'application $\alpha(.,\ell)~:A\rightarrow \mathbf{Z}\cup\{\infty\}$ est mesurable sur $A$.
\end{prop}

\begin{proof}
La seconde assertion d\'ecoule de la premi\`ere et des d\'efinitions. Le premier point est, quant \`a lui, est une cons\'equence du lemme \ref{sal et mesurable} ci-dessus qui permet de d\'ecomposer tout ensemble semi-alg\'ebrique $A$ sous la forme~:
$$
A=(M\cap A)\sqcup\left(\sqcup_i A_i\right),
$$
et du théorème 6.3.5 de \cite{s}, puisque l'assertion \ref{sal et mesurable}/3 assure que la suite $(\mu(A_i)_i$ tend vers 0.
\end{proof}

\begin{prop}
\label{sal et mesurable1}
Si $(A_\ell)_{\ell\in\mathbf{N}^r}$ est une famille semi-algébrique de $\G(X)$, alors il existe une famille semi-alg\'ebrique $(A_{\ell,i})_{(\ell,i)\in\mathbf{N}^{r+1}}$ d'ensembles stables tels que, pour tout $\ell\in\mathbf{N}^r$,
$$
\mu_{X}(A_\ell)=\sum_{i\in\mathbf{N}}\mu_{X}(A_{\ell,i})\in\widehat{\!M_k}.
$$
En particulier, si $\alpha~:\G(X)\times\mathbf{N}^{r}\rightarrow\mathbf{Z}$ est une fonction simple (cf. \S \ref{simple}), il existe une famille semi-alg\'ebrique $(A_{\ell,i})_{(\ell,i)\in\mathbf{N}^{r+1}\times\mathbf{N}}$ d'ensembles stables telle que~:
\begin{enumerate}

\item $\mu_{X}\left(\alpha(.,\ell_1,\ldots,\ell_r)^{-1}(\ell_{r+1})\right)=\sum_{i\in\mathbf{N}}\mu_{X}(A_{\ell,i}),~\forall \ell\in\mathbf{N}^{r+1}$.
 
\item $\alpha(.,\ell_1,\ldots,\ell_r)$ est constante sur $A_{\ell,i}$ pour tout $i\in\mathbf{N}$ et $\ell\in\mathbf{N}^{r+1}$.

\end{enumerate}

\end{prop}

\begin{proof}
La seconde assertion se d\'eduit de la premi\`ere en consid\'erant la famille de $\G(X)$, semi-alg\'ebrique par hypoth\`ese,
$$
A_{\ell}:=\left\{x\in\G(X)~|~\alpha(x,\ell_1,\ldots,\ell_r)=\ell_{r+1}\right\}.
$$
Pour d\'emontrer le premier point, on peut supposer que $X$ est une carte affine de $(A_\ell)_{\ell\in\mathbf{N}^r}$. On consid\`ere alors la famille $(A_{\ell,i})_{\ell\in\mathbf{N}^r}$ d\'efinie par~:
$$
A_{\ell,i}=A_\ell\cap A_i^{(\ell)},
$$
o\`u, pour chaque $\ell\in\mathbf{N}^r$, la famille semi-alg\'ebrique $(A^{(\ell)}_i)_{i\in\mathbf{N}}$ est la famille du lemme \ref{sal et mesurable} construite pour $A_\ell$. Il d\'ecoule de la construction des $A^{\ell}_i$, que la famille $(A_{\ell,i})$ est d\'efinie par une condition semi-alg\'ebrique. Le r\'esultat en d\'ecoule par stabilit\'e de la semi-alg\'ebricit\'e par intersection finie.
\end{proof}

Les r\'esultats pr\'ec\'edents se r\'esument par le diagramme commutatif suivant~:
$$
\boxed{
\xymatrix{\mathbf{B}_X^s\ar@{^(->}[r]\ar@{^(->}[d]&\mathbf{B}_X^{fs}\ar@{^(->}[r]\ar@{^(->}[d]&\mathbf{B}_X\ar@{^(->}[d]\\
\mathbf{C}_{0,\widehat{X}}\ar@{^(->}[r]&\mathbf{C}_{\widehat{X}}\ar@{^(->}[r]&\mathbf{D}_{\widehat{X}}}}
$$

\begin{rem}
En g\'en\'eral, ces notions sont diff\'erentes, \textit{i.e.} tous les cylindres (resp. ensembles mesurables) ne sont pas des ensembles faiblement stables (resp. ensembles mesurables).
\end{rem}

\subsection{Les ensembles de Presburger et les fonctions de Presburger}

On appelle \textit{ensemble de Presburger de $\mathbf{Z}^r$} un ensemble d\'efini dans $\mathbf{Z}^r$ par une combinaison bool\'eenne de conditions de la forme~:
\begin{flushleft}
$(P1)~~L(\ell_1,\ldots,\ell_r)\geq 0$\\
$(P2)~~L(\ell_1,\ldots,\ell_r)\equiv 0\mod~d$,
\end{flushleft}
pour $L(x_1,\ldots,x_r)\in\mathbf{Z}[x_1,\ldots,x_r]$ et $d\in\mathbf{N}$. Autrement dit, un ensemble de Presburger est d\'efini par une condition semi-alg\'ebrique avec $m=0$.  

\begin{exe}
\label{Pres exe}
Soit $P\subset\mathbf{Z}^r$ un ensemble de Presburger d\'efini par la condition semi-alg\'ebrique $\theta(\ell_1,\ldots\ell_r)$. Si $pr_1~:\mathbf{Z}^r\rightarrow\mathbf{Z}^{r-1}$ d\'esigne la premi\`ere projection, il d\'ecoule d'un th\'eor\`eme d'\'elimination des quantificateurs, d\^u \`a Presburger (cf. \cite{pr} ou \cite{dl1}/(2.1)), que l'ensemble $pr_1(P)$ est encore un ensemble de Presburger.
\end{exe}

Une fonction $\alpha~:\mathbf{Z}^r\rightarrow\mathbf{Z}$ est appel\'ee \textit{fonction de Presburger} si son graphe est un ensemble de Presburger de $\mathbf{Z}^{r+1}$.

\begin{exe}
Une application de la forme $\beta~:\mathbf{N}^r\rightarrow\mathbf{N}$ qui \`a $\ell\in\mathbf{N}^r$ associe l'entier $\sum_{i=1}^rb_i\ell_i +b_{r+1}$, avec $b_j\in\mathbf{N}$, est une fonction de Presburger.
\end{exe}


\section{Rationalit\'e des s\'eries de Poincar\'e motiviques}
\label{logique}

On suppose dans cette partie que $R=k[[t]]$ et que le corps $k$ est de caract\'eristique 0. Soit $X$ une $R$-vari\'et\'e plate, purement de dimension relative $d$. Nous allons montrer, principalement, un r\'esultat conjectur\'e par Du Sautoy et Loeser (cf. \cite{sl}).

\subsection{Les familles semi-alg\'ebriques born\'ees}

La notion de famille semi-alg\'ebrique born\'ee va jouer un r\^ole essentiel dans les d\'emonstrations de la proposition \ref{thm2} et du th\'eor\`eme \ref{thm3}.

\begin{defi}
\label{bornee}
\begin{itemize}
\item[1.] Un ensemble $A\subset \G(\mathbf{A}_R^N)$ est fortement born\'e s'il est d\'efinissable par une condition semi-alg\'ebrique $\theta$ et s'il existe $M\in\mathbf{N}$ tel que~:
$$
0\leq\ord_tf_i(x)\leq M,
$$
pour tout $x\in A$ et tout $i\in I$, o\`u les $f_i\in k[[t]][y_1,\ldots,y_N]$, pour $i\in I$, d\'efinissent la condition $\theta$.

\item[2.] Soit $(A_\ell)_{\ell\in\mathbf{N}^r}$ une famille semi-alg\'ebrique de $\G(X)$. On dit que la famille est born\'ee s'il existe un recouvrement fini de $X$ par des ouverts affines $U_i$ et si, pour $i\in I$, $A_\ell\cap\G(U_i)$ est d\'efinissable par une condition semi-alg\'ebrique $\theta$ et des sections r\'eguli\`eres $(h_j^{(i)})_{j\in J_i}$ de $U_i$ v\'erifiant que, via le plongement $U_i\hookrightarrow \mathbf{A}_R^{\vert J_i\vert}$ d\'efini par les $h_j^{(i)}$, l'ensemble $A_\ell\cap\G(U_i)$ est fortement born\'e, pour tout $\ell\in\mathbf{N}^r$. Autrement dit, si les polyn\^omes $f^{(i)}_j$, pour $j\in \Gamma$, d\'efinissent $\theta$, cette derni\`ere condition signifie que, pour tout $\ell\in\mathbf{N}^r$, il existe $M_\ell\in\mathbf{N}$ tel que~:
$$
0\leq\ord_tf_j^{(i)}(h_1^{(i)}(x),\ldots,h_{\vert J_i\vert}^{(i)}(x))\leq M_\ell,
$$
pour tout $x\in A_\ell\cap\G(U_i)$ et tout $j\in \Gamma$.  
\end{itemize}
\end{defi}

\begin{lem}
\label{stabilite borne}
Soit $Y$ une $R$-vari\'et\'e et $g~:Y\rightarrow X$ un $R$-morphisme de sch\'emas.
Soit $(A_\ell)_{\ell\in\mathbf{N}^r}$ une famille semi-alg\'ebrique de $\G(X)$. Alors, si $(A_\ell)_{\ell\in\mathbf{N}^r}$ est born\'ee, la famille $\left(g^{-1}(A_\ell)\right)_{\ell\in\mathbf{N}^r}$ est born\'ee.
\end{lem}

\begin{proof}
La premi\`ere assertion d\'ecoule de la proposition \ref{image inverse} et du fait que, si~:
$$
\forall y\in A_\ell, \ord_t h_j^{(i)}(y)\leq M\Rightarrow \ord_t h_j^{(i)}\circ g(x)\leq M, \forall x\in g^{-1}(A_\ell).
$$
\end{proof}

\begin{lem}
\label{atlas et borne}
Une famille $\left(A_\ell\right)_{\ell\in\mathbf{N}^r}$ semi-alg\'ebrique de $\G(X)$ est born\'ee si et seulement si elle poss\`ede un atlas semi-alg\'ebrique $(U_i)_{i\in I}$ tel que, pour tout $i\in I$, l'ensemble $A_l\cap \G(U_i)$ est fortement born\'e via un syst\`eme de coordonn\'ees sur $U_i$, pour tout $\ell\in\mathbf{N}^r$. En particulier, cette propri\'et\'e est ind\'ependante du choix d'un tel syst\`eme.
\end{lem}

\begin{proof}
La condition est clairement suffisante. Montrons qu'elle est n\'ecessaire. Soit $(U_i)_{i\in I}$ un recouvrement fini de $X$ par des ouverts affines $U_i\hookrightarrow \mathbf{A}_R^{n_i}$, vus comme sous-$R$-sch\'emas ferm\'es d'un espace affine de dimension $n_i$. Pour all\'eger les notations, on peut supposer que $X\hookrightarrow\mathbf{A}_R^n$ est \'egal \`a l'un de ces $U_i$. On peut supposer alors que la famille $(A_\ell)_{\ell\in\mathbf{N}^r}$ est d\'efinissable par $m$ sections r\'eguli\`eres $h_j$, pour $1\leq j\leq m$, et une condition semi-alg\'ebrique $\theta$. Soient $x_1,\ldots,x_n$ un syst\`eme de coordonn\'ees de $X$ et, pour tout $1\leq j\leq m$, $H_j$ le polyn\^ome tel que $h_j=H_j(x_1,\ldots,x_n)$. Comme pr\'ec\`edemment, on peut alors d\'efinir la famille $(A_\ell)_{\ell\in\mathbf{N}^r}$ \'egalement par la donn\'ee d'une condition semi-alg\'ebrique $\theta'$ et le syst\`eme de coordonn\'ees $x_1,\ldots,x_n$ (cf. \ref{atlas}). Cette donn\'ee r\'esout le probl\`eme. La derni\`ere assertion sur l'ind\'ependance du choix du syst\`eme de coordonn\'ees d\'ecoule du lemme \ref{stabilite borne} appliqu\'e \`a un morphisme de changement de coordonn\'ees.
\end{proof}

La proposition suivante prouve l'existence de pr\'esentations born\'ees pour certaines familles semi-alg\'ebriques.

\begin{prop}
\label{vers borne}
Soient $X\hookrightarrow \mathbf{A}_R^N$ une $R$-vari\'et\'e affine et $(A_\ell)_{\ell\in\mathbf{N}^r}$ une famille semi-alg\'ebrique de $\G(X)$. Supposons que, pour tout $\ell\in\mathbf{N}^r$, les $A_\ell$ sont faiblement stables. Alors cette famille peut s'\'ecrire comme une combinaison bool\'eenne de familles semi-alg\'ebriques born\'ees de $\G(X)$.
\end{prop}

\begin{proof}
Cette preuve est essentiellement une retraduction de la démonstration du lemme 2.8 de \cite{dl1}, que nous explicitons pour des commodités de lecture. On peut supposer que $X$ est une carte semi-alg\'ebrique de la famille $\left(A_\ell\right)_{\ell\in\mathbf{N}^r}$. Soient $f_1,\ldots,f_s$ les polyn\^omes intervenant dans la condition semi-alg\'ebrique d\'ecrivant la famille des $A_\ell$.

Supposons que pour $i=1,2,\ldots, e\geq s$, les fonctions $\ord_t f_i$ ne sont pas born\'ees sur $A_\nu$, pour un certain $\nu\in\mathbf{N}^r$, et que les autres fonctions $\ord_tf_{e+1},\ldots,\ord_tf_s$ sont born\'ees sur chaque $A_\ell$ par une fonction $\alpha(\ell)$. Nous allons raisonner par r\'ecurrence sur l'entier $e$. On peut supposer que $A_\ell$ est un ensemble faiblement stable de rang $\alpha(\ell)$. Supposons que $\alpha(\ell)$ est d\'esormais le plus petit entier satisfaisant les deux conditions ci-dessus. Il d\'ecoule alors des d\'efinitions et du th\'eor\`eme d'\'elimination des quantificateurs 2.1 de \cite{dl1} que l'application $\alpha~:\mathbf{N}^r\rightarrow\mathbf{N}$ est une fonction de Presburger. Un th\'eor\`eme de Schappacher (cf. \cite{sch}), utilis\'e sous la forme du th\'eor\`eme 4.3.8 de \cite{s}, assure qu'il existe une fonction de Presburger $\beta~:\mathbf{N}^r\rightarrow\mathbf{N}$, telle que, $\beta(\ell)\geq\alpha(\ell)$ pour tout $\ell\in\mathbf{N}^r$, et telle que si $x\in\G(X)$ et si~:
\begin{eqnarray}
\label{bor1}
f_1(x)\equiv f_2(x)\equiv\ldots\equiv f_e(x) \equiv 0\mod~t^{\beta(\ell)+1},
\end{eqnarray}
alors il existe $x'\in\G(X)$ tel que $x\equiv x'\mod~t^{\alpha(\ell)+1}$ et
\begin{eqnarray}
\label{bor2}
f_1(x')= f_2(x')=\ldots=f_e(x')= 0.
\end{eqnarray}
Notons que l'ensemble semi-alg\'ebrique $A_\ell$ est la r\'eunion des ensembles faiblement stables de $\G(X)$ suivants~:
$$
A_{\ell,1}:=A_\ell\cap\left\{x\in\G(X)~|~\ord_tf_1(x)\leq \beta(\ell)\right\}
$$
$$
\vdots
$$
$$
A_{\ell,e}:=A_\ell\cap\left\{x\in\G(X)~|~\ord_tf_e(x)\leq \beta(\ell)\right\}
$$
et
$$
B_{\ell}:=A_\ell\cap\left\{x\in\G(X)~|~\textrm{$x$ est solution de \ref{bor1}}\right\}.
$$
Remarquons que, gr\^ace \`a l'hypoth\`ese de r\'ecurrence, les familles $(A_{\ell,i})_{(\ell,i)}$ sont des combinaisons bool\'eennes de familles semi-alg\'ebriques born\'ees. Il ne nous reste plus qu'\`a prouver la validit\'e de l'assertion pour la famille $(B_\ell)_{\ell\in\mathbf{N}^r}$.

Consid\'erons alors la famille $(A'_\ell)_{\ell\in\mathbf{N}^r}$ d\'efinie \`a partir de $(A_\ell)_\ell$ en rempla\c{c}ant 
$$
f_1,f_2,\ldots,f_e
$$
par 0 et en ajoutant les conditions $\ord_tf_j(x)\leq \alpha(\ell)$ pour $j=e+1,\ldots,r$. Il est clair que cette famille est semi-alg\'ebrique et born\'ee et que, pour tout $\ell\in\mathbf{N}^r$, l'ensemble $A'_\ell$ est faiblement stable de rang $\alpha(\ell)$.

Si $x\in\G(X)$ satisfait \ref{bor1}, alors il existe $x'\in\G(X)$ v\'erifiant $x\equiv x'\mod~t^{\alpha(\ell)+1}$ et \ref{bor2}. Il en r\'esulte que $x\in A_\ell$ si et seulement si $x'\in A_\ell$, puisque $A_\ell$ est faiblement stable de rang $\alpha(\ell)$. Ceci est \'equivalent, par \ref{bor2}, \`a $x'\in A'_l$; et cette derni\`ere condition est v\'erifi\'ee si et seulement si $x\in A'_\ell$, puisque, pour tout $\ell\in\mathbf{N}^r$, l'ensemble $A'_\ell$ et encore faiblement stable de rang $\alpha(\ell)$. On a donc~:
$$
B_\ell=A'_l\cap\left\{x\in\G(X)~|~\textrm{$x$ est solution de \ref{bor1}}\right\}
$$
ou, ce qui revient au m\^eme,
$$
B_\ell=A'_l\backslash\bigcup_{i=1}^e\left\{x\in\G(X)~|~\ord_tf_i(x)\leq \beta(\ell)\right\}.
$$
L'assertion en d\'ecoule puisque les familles $(A'_\ell)$ et $(\left\{x\in\G(X)~|~\ord_tf_i(x)\leq \beta(\ell)\right\})$ sont born\'ees.
\end{proof}

\subsection{Un lemme de fibration pour les familles semi-alg\'ebriques born\'ees}
\label{rq2}

Soient $X$ une $R$-vari\'et\'e affine, lisse, purement de dimension relative $d$ et $\left(A_\ell\right)_{\ell\in\mathbf{N}^r}$ une famille semi-alg\'ebrique de $\G(X)$ v\'erifiant les propri\'et\'es suivantes~:

\begin{itemize}

\item[1.] $X$ est une carte affine de cette famille;

\item[2.] elle est d\'efinissable par une condition semi-alg\'ebrique $\theta$ et un syst\`eme de coordonn\'ees \textit{locales} $x_1,\ldots,x_d$, \textit{i.e.} les $d$ sections r\'eguli\`eres induisent un $R$-morphisme de sch\'emas \'etale $X\rightarrow\mathbf{A}_R^{d}$;

\item[3.] si $\theta(x_1,\ldots,x_d;\ell_1,\ldots,\ell_r)$ est d\'efini \`a partir de $s$ polyn\^omes $f_j\in R[x_1,\ldots,x_d]$, $1\leq j\leq s$, alors il existe des entiers naturels $a_{\nu j}$ v\'erifiants~:
$$
f_j(x_1,\ldots,x_d)=x_1^{a_1 j}\ldots x_d^{a_d j},
$$
pour $1\leq j\leq s$; 

\item[4.] les $tf_j$ sont encore des m\^onomes en les $x_j$;

\item[5.] pour tout $1\leq j\leq d_0$, la fonction $\ord_t x_j$ est born\'ee sur $A_\ell$ pour chaque $\ell\in \mathbf{N}^r$.
\end{itemize}

Supposons en outre que $d_0\leq d$ est le plus petit entier tel que tous les $f_j$ et les $tf_j$ s'expriment en fonction de $d_0$ coordonn\'ees, que nous supposons, quitte \`a r\'eordonner, \^etre \'egales \`a $x_1,\ldots,x_{d_0}$. Pour simplifier, supposons que $d=d_0$. La combinaison bool\'eenne des conditions du type $(SAL3)$ de $\theta$ d\'efinit un ensemble constructible de $\mathbf{G}_{m,k}^{d_0}$, que l'on note $C\subset\mathbf{G}_{m,k}^{d_0}$. Comme dans le paragraphe \S3, on note $D_j$ l'hypersurface de $X$~:
$$
D_j=V\left(f_j(x_1,\ldots,x_{d_0})\right)\hookrightarrow X.
$$

\begin{prop}
\label{pi0 fini}
Sous les hypoth\`eses ci-dessus, il existe une famille finie $(\mathsf{C_q})_{q\in Q}$ d'\'el\'ements de $\!M_k$, une fonction de Presburger $\alpha~:\mathbf{Z}^{d_0}\rightarrow\mathbf{Z}$ et, pour tout $q\in Q$, un ensemble de Presburger $P_q$ tels que, pour tout $\ell\in\mathbf{N}^r$,
$$
\mu_\!X(A_\ell)=\sum_{q\in Q}[\mathsf{C_q}]\sum_{n\in\mathbf{N}^{d_0}, n\in P_q(n,\ell)}\mathbf{L}^{-\alpha(n)}.
$$
\end{prop}

\begin{proof}
Soit $Q$ l'ensemble des parties de l'ensemble $\{1,\ldots,d_0\}$. Soit $\alpha~:\mathbf{Z}^{d_0}\rightarrow \mathbf{Z}$ l'application qui \`a $n\in\mathbf{Z}^{d_0} $ associe l'entier $d+\sum_{i=1}^{d_0} n_i$. Cette application est de Presburger par d\'efinition. Si $J\in Q$, on pose dans $\!M_k$~:
$$
\mathsf{C}_{J}:=[\mathsf{D}_J^\circ\times_{\mathbf{G}_{m,k}^{d_0-\vert J\vert}} C].
$$

La condition $\theta'(n_1,\ldots,n_{d_0},\ell_1,\ldots,\ell_r)$ d\'efinie par,
$$
\exists x\in\theta,~~n_i=\ord_tx_i(\varphi_x)~\forall 1\leq i\leq d_0,
$$
est une condition semi-alg\'ebrique gr\^ace au th\'eor\`eme d'\'elimination des quantificateurs de Pas (cf. th\'eor\`eme de Pas 2.1 de \cite{dl1}). Cette condition d\'efinit donc un ensemble de Presburger de $\mathbf{Z}^{r+\vert J\vert}$. Si $J$ est le support de $n$, on le note $P_J(n,\ell)$.

La section r\'eguli\`ere au-dessus de $X$ $x_i$ d\'efinit un $R$-morphisme que l'on note encore $x_i~:X\rightarrow\mathbf{A}^1_R$. Pour $n\in\mathbf{N}^{d_0}$, on peut alors d\'efinir une application~:
$$
ac_n~:\cap_{i=1}^{d_0}\ord_t x_i^{-1}(n_i)\rightarrow\mathbf{G}^{d_0}_{m,k},
$$
qui consiste \`a associer \`a $x\in\G(X)$ le $d_0$-uplet, dont la $i$-\`eme coordonn\'ee est le premier coefficient non nul de $x_i(\varphi_x)$, vu comme \'el\'ement de $\kappa(x)[[t]]$. Si $J$ est le support de $n$, posons
$$
\!W_n:=\pi_{0,X}^{-1}(\mathsf{D}_J^\circ)\cap ac_n^{-1}(C).
$$
Il est clair que cet ensemble est un $(\vert n\vert:=\sum_{i=1}^{d_0} n_i)$-cylindre de $\G(X)$. Le choix des $x_i$ permet d'identifier $\G_{\vert n\vert}(\!X)$ et $\mathbf{A}^{(\vert n\vert +1)d}_k\times_{\mathbf{A}_k^d} X_0$. Le $k$-morphisme de sch\'emas $\theta^n_0~:\G_{\vert n\vert}(X)\rightarrow\G_0(X)$ induit une application~:
$$
\psi_n~:\pi_{\vert n\vert,X}(\!W_n)\rightarrow \mathsf{D}_J^\circ\times_{\mathbf{G}_{m,k}^{d_0-\vert J\vert}} C.
$$
Comme les $f_j$ sont des mon\^omes en les $x_i$, il est facile de voir, comme dans la preuve du lemme \ref{calcul dcn} que, pour tout $x\in (\mathsf{D}_J^\circ\times_{\mathbf{G}_{m,k}^{d_0-\vert J\vert}} C)(\kappa(x))$, les fibres de cette application sont isomorphes \`a $\mathbf{A}_{\kappa(x)}^{\vert n\vert(d-1)}$. Il d\'ecoule donc du th\'eor\`eme 4.2.3 de \cite{s} que cette application est une fibration localement triviale. En particulier, on a dans $\!M_k$ la relation~:
$$
[\pi_{\vert n\vert,X}(\!W_n)]=\mathsf{C}_{J}\times\mathbf{L}^{\vert n\vert(d-1)}.
$$

Il nous reste \`a calculer le volume $\mu_X(A_\ell)$. Par d\'efinition, 
$$
A_\ell=\bigsqcup_{q\in Q}\bigsqcup_{m\in\mathbf{N}^{\vert J\vert},P_q(m,\ell)}\!W_n.
$$
Par additivit\'e de la mesure et avec les notations du d\'ebut,
$$
\mu_X(A_\ell)=\sum_{q\in Q}\mathsf{C}_q\sum_{m\in\mathbf{N}^{\vert J\vert},P_q(m,\ell)}\mathbf{L}^{-\alpha(m)}.
$$
\end{proof}

\subsection{L'\'enonc\'e des r\'esultats}

Le r\'esulat le plus g\'en\'eral de ce paragraphe est le th\'eor\`eme \ref{thm3}. Il a \'et\'e propos\'e par Du Sautoy et Loeser. Commen\c{c}ons par rappeler un r\'esultat d\'emontr\'e par Denef et Loeser~:

\begin{lem}
\label{dl1}
Soit $P$ un ensemble de Presburger de $\mathbf{Z}^m$ et $\varphi_i~:\mathbf{Z}^m\rightarrow\mathbf{N}$, pour $1\leq i\leq m$, des fonctions de Presburger. Supposons que les fibres de l'application $\varphi~: P\rightarrow\mathbf{N}^r$ d\'efinie par $i\mapsto(\varphi_1(i),\ldots,\varphi_r(i))$ sont finies. Alors la s\'erie $f(X):=\sum_{i\in P}X^{\varphi(i)}$, avec $X=(X_1,\ldots,X_r)$, appartient au sous-anneau de $\mathbf{Z}[[X]]$ engendr\'e par $\mathbf{Z}[X]$ et les s\'eries $(1-X^c)^{-1}$, avec $c\in\mathbf{N}^r\backslash\{0\}$.
\end{lem}

\begin{proof}
C'est le lemme 5.2 de \cite{dl1}.
\end{proof}

Cet \'enonc\'e va permettre de d\'emontrer la proposition suivante~:

\begin{prop}
\label{pro2}
Soit $X$ une $R$-vari\'et\'e plate, purement de dimension relative $d$. Soit $A_n$, $n\in\mathbf{Z}^r$, une famille semi-alg\'ebrique de $\G(X)$. Soit $\alpha~:\G(X)\times\mathbf{Z}\rightarrow\mathbf{N}$ une fonction simple. Supposons que $A_n\cap\G(X_{\mathrm{sing}})=\emptyset$ et que, pour tout $n\in\mathbf{N}^r$, les ensembles $A_n\cap (\alpha(.,n)^{-1}(m))$ sont stables quel que soit $m\in\mathbf{N}$. Alors la  s\'erie de $\overline{\!M_k}[[T]]$
\begin{eqnarray}
\sum_{n\in\mathbf{N}^r}\int_{A_n}\mathbf{L}^{-\alpha(.,n)}d\mu_{X}T^n,
\end{eqnarray}
en la variable $T=(T_1,\ldots,T_r)$, appartient au sous-anneau de $\overline{\!M_k}[[T]]$ engendr\'e par $\overline{\!M_k}[T]$, par les s\'eries $(\mathbf{1}-\mathbf{L}^{-a}T^b)^{-1}$ et $(\mathbf{L}^{i}-\mathbf{1})^{-1}$, avec $a\in\mathbf{N}$ et $b\in\mathbf{N}^r\backslash\{0\}$ et $i\in\mathbf{N}$.
\end{prop}

\begin{proof}
\label{dem thm2}

Si $n\in\mathbf{N}^r$ et $m\in\mathbf{N}$, on pose~:
$$
A_{n,m}:=\left\{x\in A_n~|~\alpha(x,n)=m\right\}.
$$
Par d\'efinition de l'int\'egrale,
\begin{eqnarray}
\label{fin1}
\sum_{n\in\mathbf{N}^r}\left(\int_{A_n}\mathbf{L}^{-\alpha(.,n)}d\mu_X\right)T^n=\sum_{n\in\mathbf{N}^r}\left(\sum_{m\in\mathbf{N}}\mu_{X}(A_{n,m})\mathbf{L}^{-m}\right)T^n.
\end{eqnarray}
Comme, pour tout $n\in\mathbf{N}^r$, l'application $\alpha(.,n)$ est \`a valeurs dans $\mathbf{N}$, remarquons que~:
$$
A_n\subset\bigsqcup_{m\in\mathbf{N}} A_{n,m},
$$
pour tout $n\in\mathbf{N}^r$. En particulier, comme $A_n$ est faiblement stable, on peut supposer, par quasi-compacit\'e de la topologie constructible (cf. lemme 4.3.10 de \cite{s}), que $A_n$ est recouvert par un nombre fini de $A_{n,m}$. Ceci justifie que la s\'erie (\`a support fini!)~:
$$
\sum_{m\in\mathbf{N}}\mu_{X}(A_{n,m})\mathbf{L}^{-m},
$$
converge dans $\widehat{\!M_k}$ (et m\^eme dans $\overline{\!M_k}$).

\smallskip

Le th\'eor\`eme de changement de variables appliqu\'e \`a une r\'esolution de N\'eron plong\'ee de $(X,X_{\mathrm{sing}})$ permet de se ramener au cas o\`u $X\rightarrow R$ est lisse. L'additivit\'e de la mesure $\mu_X$ (cf. lemme 7.0.7 \cite{s}) permet de supposer que $X$ est une $R$-vari\'et\'e affine. En outre, gr\^ace au lemme \ref{vers borne}, les relations d'additivit\'e et la stabilit\'e de la propri\'et\'e d'\^etre born\'e par intersection finie permettent de supposer que la famille semi-alg\'ebrique de $\G(X)$~:
$$
\left(A_{n,m}\right)_{(n,m)\in\mathbf{N}^r\times\mathbf{N}}
$$
est born\'ee. Enfin, toujours par additivit\'e de la mesure et gr\^ace au lemme \ref{atlas et borne}, on peut se ramener au cas o\`u $X$ est une carte semi-alg\'ebrique  pour la famille des $A_{n,m}$.

\smallskip

Si $F$ est le produit de tous les polyn\^omes $f_i$,  $1\leq i\leq m$, qui interviennent dans la description des conditions du type $(SAL1)$, $(SAL2)$ et $(SAL3)$ d\'efinissant la famille $(A_{m,n})_{(n,m)}$, on peut consid\'erer $h~:Y\rightarrow X$ une r\'esolution de N\'eron plong\'ee de $(X,V(tF))$. Le th\'eor\`eme de changement de variables donne alors la formule~:
$$
\mu_X(A_{n,m})=\sum_{e\in\mathbf{N}}\mu_{Y}\left(h^{-1}(A_{n,m})\cap\ord_t(\mathrm{Jac})_h^{-1}(e)\right)\mathbf{L}^{-e}.
$$
Comme pr\'ec\'edemment, on peut supposer, puisque la topologie constructible est quasi-compacte, que cette somme est \`a support fini, donc convergente dans $\overline{\!M_k}$. Par construction, il existe un recouvrement fini de $Y$ par des ouverts affines $(V_i)_{i\in I}$ et, pour tout $i\in I$, $d$ sections $z_1,\ldots,z_d$ sur $V_i$, induisant un $R$-morphisme \'etale de sch\'emas $V_i\rightarrow\mathbf{A}^d_R$, et telles que, pour tout $1\leq j\leq m$, $f_j$ s'exprime comme un mon\^ome en les $z_1,\ldots,z_d$. En outre, on peut supposer que, sur chaque $V_i$, la famille semi-alg\'ebrique des $h^{-1}(A_{n,m})$ est d\'efinissable par une condition semi-alg\'ebrique $\theta$ et les $z_j$ et, qu'ainsi d\'efinie, cette famille est born\'ee, puisque celle des $A_{n,m}$ l'est. Par ailleurs, comme $X$ est lisse sur $R$, la famille semi-alg\'ebrique de $\G(Y)$,
$$
\left(\ord_t(\mathrm{Jac})_h^{-1}(e)\right)_{e\in\mathbf{N}}
$$
est d\'efinissable sur chaque $V_i$ par une section r\'eguli\`ere $g(z_1,\ldots,z_d)=z_1^{b_1}\ldots z_d^{b_d}$, avec $b_j\geq 0$. En particulier, cette famille est born\'ee. Supposons que $g$ et les $f_j$ sont des mon\^omes en $z_1,\ldots,z_{d_0}$ et que l'entier $d_0\leq d$ est le plus petit entier qui v\'erifie cette condition. On peut supposer, par additivit\'e de la mesure $\mu_Y$, que $Y$ est \'egal \`a l'un des $V_i$.

\smallskip

Il d\'ecoule de l'\'etude men\'ee au \S \ref{rq2} et de la proposition \ref{pi0 fini} qu'il existe une famille finie d'\'el\'ements $(\mathsf{C}_q)_{q\in Q}$ et un ensemble de Presburger $P$ de $\mathbf{Z}^{r+1+d_0}$ tel que, pour tout $n\in\mathbf{N}^r$ et tout $m\in \mathbf{N}$,
\begin{eqnarray}
\label{fin2}
\mu_{X}(A_{n,m})=\sum_{q\in Q}\mathsf{C}_q\sum_{P(\ell,n,m)}\mathbf{L}^{-d-\sum_{i=1}^{d_0}\ell_i}\mathbf{L}^{-\beta(\ell)},
\end{eqnarray}
L'application $\beta$ est la forme lin\'eaire \`a coefficients dans $\mathbf{N}$, qui, \`a $\ell\in\mathbf{N}^{d_0}$, associe l'entier naturel $\sum_{i=1}^{d_0} b_i\ell_i$. 
Les \'egalit\'es de \ref{fin1} et \ref{fin2} entrainent donc que la s\'erie~:
$$
\sum_{n\in\mathbf{N}^r}\int_{A_n}\mathbf{L}^{-\alpha(.,n)}d\mu_{\widehat{X}}T^n
$$
est une $\!M_k$-combinaison lin\'eaire de s\'eries formelles de la forme $f(\mathbf{L}^{-1},T_1,\ldots,T_r)$, avec $f(U,X_1,\ldots,X_r)\in\mathbf{Z}[[U,X_1,\ldots,X_r]]$ de la forme de celle du lemme \label{dl1} ci-dessus.
\end{proof}

Les hypoth\`eses de de la proposition \ref{pro2} permettent en r\'ealit\'e d'affiner l'\'enonc\'e en~:

\begin{thm}
\label{thm2}
Soit $X$ une $R$-vari\'et\'e plate, purement de dimension relative $d$. Soit $A_n$, $n\in\mathbf{Z}^r$, une famille semi-alg\'ebrique de $\G(X)$. Soit $\alpha~:\G(X)\times\mathbf{Z}\rightarrow\mathbf{N}$ une fonction simple. Supposons que $A_n\cap\G(X_{\mathrm{sing}})=\emptyset$ et que, pour tout $n\in\mathbf{N}^r$, les ensembles $A_n\cap (\alpha(.,n)^{-1}(m))$ sont stables quel que soit $m\in\mathbf{N}$. Alors la  s\'erie de ${\!M_k}[[T]]$
\begin{eqnarray}
\sum_{n\in\mathbf{N}^r}\int_{A_n}\mathbf{L}^{-\alpha(.,n)}d\mu_{X}T^n,
\end{eqnarray}
en la variable $T=(T_1,\ldots,T_r)$, appartient au sous-anneau de ${\!M_k}[[T]]$ engendr\'e par ${\!M_k}[T]$, par les s\'eries $(\mathbf{1}-\mathbf{L}^{-a}T^b)^{-1}$ et $(\mathbf{L}^{i}-\mathbf{1})^{-1}$, avec $a\in\mathbf{N}$ et $b\in\mathbf{N}^r\backslash\{0\}$ et $i\in\mathbf{N}$.
\end{thm}

\begin{proof}
Commen\c{c}ons par donner un sens \`a l'\'enonc\'e de ce th\'eor\`eme. En effet, \textit{a priori}, 
$$
\mu_X(A_n)\in\overline{\!M_k}.
$$
Toutefois, si $A$ est un $\nu$-cylindre contenu dans $\G(X)\backslash\G(X_{\mathrm{sing}})$, il d\'ecoule des d\'efinitions (cf. \cite{s}) que~:
$$
\mu_X(A)=\iota(\mu_{0,X}(A)),
$$
o\`u $\iota~:\!M_k\rightarrow\widehat{\!M_k}$ est le morphisme d'anneaux canonique et $\mu_{0,X}$ l'application additive sur $\mathbf{C}_{0,\widehat{X}}$ d\'efinie par~:
$$
\mu_{0,X}(A)=\left[\pi_{\nu,X}(A)\right]\mathbf{L}^{-(\nu+1)d}.
$$
En identifiant, $\mu_X$ et $\mu_{0,X}$, on peut donc supposer que~:
$$
\mu_X(A_n)\in{\!M_k}.
$$
Passons \`a la d\'emonstration proprement dite. Remarquons que si $h~:Y\rightarrow X$ est une r\'esolution de N\'eron plong\'ee de $(X,Z)$, $Z$ une sous-$R$-vari\'et\'e ferm\'ee de $X$ contenant $X_{\mathrm{sing}}$, la condition, v\'erifi\'ee par $h$, que~:
$$
h^{-1}(Z)\supset \ord_t(\mathrm{Jac})_h^{-1}(\infty)
$$ 
entra\^\i ne en particulier que, pour tout cylindre $A\subset\G(X)\backslash\G(Z)$, il existe un nombre fini d'entiers naturels $e\in I$ et un entier naturel $e'$ tels que~:
$$
h^{-1}(A)\subset\left(\bigsqcup_{e\in I}\ord_t(\mathrm{Jac})_h^{-1}(e)\right)\cap h^{-1}(\G^{e'}(X)).
$$
Il d\'ecoule alors de la premi\`ere remarque ci-dessus et du lemme 8.1.3 de \cite{s} que~:
$$
\mu_X(h^{-1}(A))\in\!M_k.
$$
Cette derni\`ere remarque permet d'interpr\'eter le lemme 8.1.3 de \cite{s} comme un th\'eor\`eme de changement de variables dans $\!M_k$.

La preuve du th\'eor\`eme \ref{thm2} devient alors une simple r\'eadaptation de la d\'emonstration de la proposition \ref{pro2} (cf. \S \ref{dem thm2}) en substituant, \`a l'usage du th\'eor\`eme de changement de variable `` classique '', celui du lemme 8.1.3 de \cite{s}. 
\end{proof}

\begin{cor}
Si $k$ est un corps de caract\'eristique 0 et $X$ une $k[[t]]$-vari\'et\'e plate, purement de dimension relative $d$. Alors 
$$
\mu_{X}(\G(X))\in\overline{\!M_k}\left[\left(\frac{\mathbf{L}-\mathbf{1}}{\mathbf{L}^{i}-\mathbf{1}}\right)_{i\geq 1}\right].
$$
\end{cor}

\begin{proof}
Soit $h~:Y\rightarrow X$ une r\'esolution de N\'eron plong\'ee de $X$. Le th\'eor\`eme de changement de variables assure alors que~:
$$
\mu_X(\G(X))=\int_{\G(Y)}\mathbf{L}^{-\ord_t(\mathrm{Jac})}d\mu_Y:=\mathbf{L}^{-d}\sum_{n\in\mathbf{N}}\mu_Y(\ord_t(\mathrm{Jac})^{-1}(n))\mathbf{L}^{-nd}.
$$
On applique la proposition \ref{pro2} ci-dessus \`a la famille~:
$$
A_n:=\left(\ord_t(\mathrm{Jac})^{-1}(n)\right)_{n\in\mathbf{N}}
$$
et \`a la s\'erie~:
$$
V(T):=\sum_{n\in\mathbf{N}} T^n\int_{A_n}\mathbf{L}^{-0}d\mu_Y.
$$
Le r\'esultat en d\'ecoule du fait que, gr\^ace \`a la premi\`ere \'egalit\'e ci-dessus, $V(\mathbf{L}^{-1})=\mu_X(\G(X))$.

\end{proof}

\begin{cor}
\label{coro poincare}
Soit $Y$ une $R$-vari\'et\'e. Soit $A$ un ensemble semi-alg\'ebrique de $\G(Y)$. Alors la s\'erie~:
\begin{eqnarray}
P_A(T):=\sum_{n=0}^{+\infty}[\pi_{n,\widehat{Y}}(A)]T^n ~\in~\!M_k[[T]],
\end{eqnarray}
appartient au sous-anneau de ${\!M_k}[[T]]$ engendr\'e par ${\!M_k}[T]$, les s\'eries $(\mathbf{1}-\mathbf{L}^{-a}T^b)^{-1}$ et $(\mathbf{L}^{i}-\mathbf{1})^{-1}$, avec $a\in\mathbf{N}$ et $b\in\mathbf{N}\backslash\{0\}$ et $i\in\mathbf{N}$.
\end{cor}

\begin{proof}
On peut supposer que $Y$ est une sous-$R$-vari\'et\'e d'une vari\'et\'e $X$, lisse, connexe, de dimension $d$ et de fibre sp\'eciale r\'eduite. Par d\'efinition,
$$
\int_{\pi_{n,X}^{-1}(\pi_{n,Y}(A))}\mathbf{L}^{-0}d\mu_{X}=[\pi_{n,Y}(A)]\mathbf{L}^{-(n+1)d}.
$$
La s\'erie $\mathbf{L}^{-d}\sum_{n=0}[\pi_{n,X}(A)](\mathbf{L}^{-d}T)^n$ appartient (cf. lemme \ref{Poincare} ci-dessous et th\'eor\`eme \ref{thm2}) au sous-anneau de ${\!M_k}[[T]]$ engendr\'e par ${\!M_k}[T]$ et les s\'eries $(\mathbf{1}-\mathbf{L}^{-a}T^b)^{-1}$, avec $a\in\mathbf{N}$ et $b\in\mathbf{N}\backslash\{0\}$.
\end{proof}

Le lemme ci-dessous est une g\'en\'eralisation de l'exemple \ref{ordjac}/1~:

\begin{lem}
\label{Poincare}
Soit $X$ une $R$-vari\'et\'e lisse. Soient $Y\hookrightarrow X$ une sous-$R$-vari\'et\'e localement ferm\'ee de $X$ et $A\subset\G(\widehat{Y})$ un ensemble semi-alg\'ebrique de $\G(\widehat{Y})$. Alors la famille~:
$$
\Big(\pi_{n,\widehat{X}}^{-1}\big(\pi_{n,\widehat{Y}}(A)\big)\Big)_{n\in\mathbf{N}}
$$
est une famille semi-alg\'ebrique d'ensembles stables de $\G(X)$.
\end{lem}

\begin{proof}
La stabilit\'e des ensembles 
$$
B_n:=\pi_{n,\widehat{X}}^{-1}\big(\pi_{n,\widehat{Y}}(A)\big)
$$
d\'ecoule de la lissit\'e de $X$. Il existe $U\hookrightarrow Y$ un ouvert affine de $X$ tel que $U\cap Y$ est une carte semi-alg\'ebrique de $A$ et tel que $Y\cap U\hookrightarrow U$ est une immersion ferm\'ee. On peut donc supposer que $X\hookrightarrow\mathbf{A}_R^N$ est affine, que $Y$ est une carte affine de $A$ et que $Y$ est une sous-$R$-vari\'et\'e ferm\'ee de $X$. Soit $\theta$ la condition semi-alg\'ebrique d\'efinissanr la $A$. Consid\'erons $\theta'$ la condition semi-alg\'ebrique obtenue en ajoutant \`a $\theta$ les conditions~:
$$
\ord_t f_j(x_j,y_j)\geq L(n),
$$
avec $f_j(x_j,y_j)=x_j-y_j$, pour $1\leq j\leq N$, et $L(x)=x+1$, et~:
$$
\ord_t g_i=\infty,
$$
o\`u les $g_i$, pour $1\leq i\leq m$, sont une pr\'esentation de sections r\'eguli\`eres de $X$ d\'efinissant $Y$ comme sous-vari\'et\'e ferm\'ee. Le th\'eor\`eme d'\'elimination des quantificateurs 2.1 de \cite{dl1} implique, comme pour l'exemple \ref{ordjac}/1, que la famille des $B_n$ est une semi-alg\'ebrique, d\'efinissable \`a partir de la condition $\theta'$. 
\end{proof}

\begin{cor}
Soit $X$ une $R$-vari\'et\'e. La s\'erie de Poincar\'e motivique associ\'ee \`a $X$,
$$
P_X(T):=\sum_{n\in\mathbf{N}}[\pi_{n,X}(\G(X))]T^n\in\!M_k[[T]]
$$
appartient au sous-anneau de $\!M_k[[T]]$ engendr\'e par $\!M_k[T]$, les s\'eries $(\mathbf{1}-\mathbf{L}^{-a}T^b)^{-1}$ et $(\mathbf{L}^{i}-\mathbf{1})^{-1}$ , avec $a\in\mathbf{N}$ et $b\in\mathbf{N}\backslash\{0\}$ et $i\in\mathbf{N}$.
\end{cor}

\begin{thm}
\label{thm3}
Soit $X$ une $R$-vari\'et\'e plate, purement de dimension relative $d$. Soit $A_n$, $n\in\mathbf{Z}^r$, une famille semi-alg\'ebrique de $\G(X)$ et soit $\alpha~: \G(X)\times\mathbf{Z}^r\rightarrow\mathbf{N}$ une fonction simple. Alors, si, pour tout $n\in\mathbf{N}^r$, l'application $\alpha(.,n)~:A_n\rightarrow\mathbf{N}$ est exponentiellement int\'egrable, la s\'erie 
\begin{eqnarray}
\label{eheh}
\sum_{n\in\mathbf{N}^r}\int_{A_n}\mathbf{L}^{-\alpha(.,n)}d\mu_{X}T^n,
\end{eqnarray}
en $T=(T_1,\ldots,T_r)$ appartient au sous-anneau de $\widehat{\!M_k}[[T]]$ engendr\'e par $\overline{\!M_k}[T]$, par les s\'eries $(\mathbf{L}^i-\mathbf{1})^{-1}$ et $(\mathbf{1}-\mathbf{L}^{-a}T^b)^{-1}$, o\`u $i\in\mathbf{N}\backslash\{0\}$, $a\in\mathbf{N}$ et $b\in\mathbf{N}^r\backslash\{0\}$.
\end{thm}

\begin{proof}
La d\'emonstration de ce th\'eor\`eme est analogue \`a celle de la proposition \ref{pro2}. L'unique diff\'erence vient du fait qu'ici il n'est plus n\'ecessaire de consid\'erer des familles born\'ees. Le r\'esultat d\'ecoule donc du lemme \ref{dl1} et de la proposition \ref{pi0 fini}.
\end{proof}

\section{La rationalit\'e des fonctions Z\^eta d'Igusa motiviques}
\label{3}

Nous supposons que $R=k[[t]]$, avec $k$ un corps de caract\'eristique 0. 

\subsection{Les diviseurs \`a croisements normaux}

\subsubsection{}
\label{ici}
Soit $Y$ un sch\'ema. Soit $D\hookrightarrow Y$ un diviseur de $Y$. On dit que $D$ est un diviseur \`a croisements normaux, si~:
\begin{enumerate}

\item[\textsl{(i)}] pour tout $s\in D$, l'anneau local $\!O_{Y,s}$ est r\'egulier;

\item[\textsl{(ii)}] pour tout $s\in D$, il existe un syst\`eme r\'egulier de param\`etres en $s$, $t_1,\ldots,t_d$, tel qu'une \'equation de $D$, au voisinage de $s$, est $t_1^{n_1}\ldots t_d^{n_d}=0$, avec $n_i\geq 0$ pour tout $1\leq i\leq d$.

\end{enumerate}

Si, en outre, les $D_i$, pour tout $i\in I$, désignent les composantes irréductibles de $D$ et si les sch\'emas $D_i$ sont r\'eguliers, on dit que le diviseur $D$ est \`a croisements normaux stricts dans $Y$.

\begin{lem}
\label{dcn}
Soit $Y$ une $R$-vari\'et\'e lisse. Soit $D\hookrightarrow Y$ un diviseur \`a croisements normaux (resp. stricts) de $Y$. Si $E\hookrightarrow Y$ est un diviseur de $Y$ contenu dans $D$, alors $E$ est à croisements normaux (resp. stricts).
\end{lem}

\begin{proof}
L'assertion d\'ecoule directement de la définition \ref{ici}.
\end{proof}

Si $\mathsf{Z}\hookrightarrow \mathsf{Y}$ est une sous-$k$-vari\'et\'e ferm\'ee de $\mathsf{Y}$, on note $(\mathsf{Z}_j)_{1\leq j\leq r}$ l'ensemble des composantes irr\'eductibles de $\mathsf{Z}$. Pour $I\subset\{1,\ldots,r\}$, on note $\mathsf{Z}_I=\cap_{i\in I} \mathsf{Z}_i$, si $I\not=\emptyset$, et $\mathsf{Y}$ sinon. On pose $\mathsf{Z}^\circ_I:= \mathsf{Z}_I\backslash \cup_{i\in\{1,\ldots,r\}\backslash J} \mathsf{Z}_i$. En particulier, $\mathsf{Z}^\circ_\emptyset=\mathsf{Y}\backslash \mathsf{Z}$. Il d\'ecoule \'egalement de la d\'efinition que $\mathsf{Y}=\sqcup_{I\subset\{1,\ldots,r\}} \mathsf{Z}^\circ_I$.

\subsection{Les notations}

Dans ce paragraphe, nous adoptons les notations et conventions suivantes~:

\begin{itemize}

\item[1.] Si $\mathsf{Z}\hookrightarrow \mathsf{Y}$ est une sous-$k$-vari\'et\'e ferm\'ee de $\mathsf{Y}$, on note $(\mathsf{Z}_j)_{1\leq j\leq r}$ l'ensemble des composantes irr\'eductibles de $\mathsf{Z}$. Pour $I\subset\{1,\ldots,r\}$, on note $\mathsf{Z}_I=\cap_{i\in I} \mathsf{Z}_i$, si $I\not=\emptyset$, et $\mathsf{Y}$ sinon. On pose $\mathsf{Z}^\circ_I:= \mathsf{Z}_I\backslash \cup_{i\in\{1,\ldots,r\}\backslash J} \mathsf{Z}_i$. En particulier, $\mathsf{Z}^\circ_\emptyset=\mathsf{Y}\backslash \mathsf{Z}$. Il d\'ecoule \'egalement de la d\'efinition que $\mathsf{Y}=\sqcup_{I\subset\{1,\ldots,r\}} \mathsf{Z}^\circ_I$.

\item[2.] Si $Y$ est une $R$-vari\'et\'e alg\'ebrique lisse et $D\hookrightarrow Y$ un diviseur de $Y$, on note $r$ le nombre de composantes connexes de $\mathsf{Y}$. On suppose que $D$ poss\`ede $m$ composantes irr\'eductibles que l'on note $D_i$, pour $1\leq i\leq m$. On pose $I$ l'ensemble des $i$ tels que $D_i\subset \mathsf{Y}$. Si $x\in D_i$, on note $I_x$ l'ensemble des $i$ tels que $x\in D^\circ_{I_x}$. 

\item[3.] Si $n\in\mathbf{N}^m$, on note $n_i$ ses coordonn\'ees, pour $1\leq i\leq m$, et $J$ son support, \textit{i.e.} l'ensemble des $i$ tels que $n_i\not=0$. Si $D$ est un diviseur de $Y$, on dit que $n$ ne v\'erifie pas la \textit{propri\'et\'e $\mathfrak{P}_D$}, ce que l'on note $n\not\in\mathfrak{P}_D$, si 
$$
\left\{\begin{array}{l}
\textrm{$r\geq \vert I\vert\geq 1$ et si $\exists i\in I$ tel que  $n_{i}\not\in\{0,1\}$;}\\
\textrm{ou}\\
\textrm{$r\geq \vert I\vert\geq 1$ et si $\exists i,j\in I$ tels que $n_{i}=n_{j}=1$;}\\ 
\textrm{ou}\\
\textrm{$\vert I\vert=r$ et si $\forall i\in I, n_{i}\not=1$}.
\end{array}\right.
$$
Dans le cas contraire, on dit que $n$ v\'erifie la propri\'et\'e $\mathfrak{P}_D$, ce que l'on note $n\in\mathfrak{P}_D$.

\item[4.] nous ne nous int\'eressons qu'aux diviseurs $D$ tels que $\widehat{D}\not=\emptyset$. 

\end{itemize}

\subsection{La fonction Z\^eta motivique associ\'ee \`a un diviseur \`a croisements normaux stricts d'une vari\'et\'e lisse}

Soit $Y$ une $R$-vari\'et\'e lisse, purement de dimension relative $d$. Soit $D\hookrightarrow Y$ un diviseur \`a croisements normaux stricts, réduit. Supposons en outre que $tD$ est encore un diviseur à croisements normaux stricts dans $Y$. La fonction Z\^eta (d'Igusa) motivique associ\'ee \`a $D$ est la s\'erie de ${\!M_k}[[T_1,\ldots,T_m]]$ d\'efinie par~:
\begin{eqnarray}
Z_D(T):=\sum_{n:=(n_1,\ldots,n_m)}[\pi_{|n|,Y}\big(\cap_{i=1}^m\mathrm{mult}_{D_i}^{-1}(n_i)\big)]\mathbf{L}^{-|n|d}T^n,
\end{eqnarray}
avec les notations $T^n:=T_1^{n_1}\ldots T_m^{n_m}$ et $\vert n\vert=\sum_{i=1}^m n_i$ pour tout $n\in\mathbf{N}^m$. Si $x\in\mathsf{D}(k)$, on associe au couple $(D,x)$ la fonction Z\^eta de \textit{support $x$}~:
$$
Z_{D,x}(T):=\sum_{n:=(n_1,\ldots,n_m)}\Big[\pi_{|n|,Y}\Big(\big(\cap_{i=1}^m\mathrm{mult}_{D_i}^{-1}(n_i)\big)\cap \pi_{0,Y}^{-1}(x)\Big)\Big]\mathbf{L}^{-|n|d}T^n.
$$

Commen\c{c}ons par d\'emontrer un lemme g\'en\'eral qui interviendra dans les calculs des diff\'erents r\'esultats ci-dessous~:

\begin{lem}
\label{calcul serie}
Soient $A$ un anneau commutatif et $m$ un entier naturel non nul. Si $a\in A$ est un \'el\'ement inversible, alors la s\'erie~:
$$
S_a(T):=\sum_{n\in(\mathbf{N}\backslash\{0\})^m} a^{-\vert n\vert} T^n=\prod_{i=1}^m\frac{a^{-1}T_i}{1-a^{-1}T_i}~\in~A[[T_1,\ldots,T_m]].
$$
\end{lem}

\begin{proof}
Il d\'ecoule de la d\'efinition que~:
$$
S_a(T)=\sum_{n_1>0,\ldots,n_m>0}\left(\prod_{i=1}^{m}a^{-n_i}T_i^{n_i}\right).
$$
Par suite,
$$
S_a(T)=\prod_{i=1}^m\left(\sum_{n_i>0}a^{-n_i}T_i^{n_i}\right),
$$
ou ce qui revient au m\^eme
$$
S_a(T)=\prod_{i=1}^m\left(\frac{1}{1-a^{-1}T_i}-1\right),
$$
soit~:
$$
S_a(T)=\prod_{i=1}^m\frac{a^{-1}T_i}{1-a^{-1}T_i},
$$
\end{proof}

\begin{lem}
\label{calcul dcn}
Soit $n\in\mathbf{N}^m$. Alors on a, dans $\!M_k$, la relation
$$
[\pi_{|n|,Y}\big(\cap_{i=1}^m\mathrm{mult}_{D_i}^{-1}(n_i)\big)]=\lbrack\mathsf{D}^{\circ}_J\rbrack\times(\mathbf{L}-\mathbf{1})^{\vert J\vert}\times \mathbf{L}^{\vert n\vert (d-1)}, 
$$
si $n\in\mathfrak{P}_D$. Sinon,
$$
[\pi_{|n|,Y}\big(\cap_{i=1}^m\mathrm{mult}_{D_i}^{-1}(n_i)\big)]=0.
$$

\end{lem}

\begin{proof}

Soit $n\in\mathbf{N}^m$. Posons~:
$$
A_n:=\cap_{i=1}^m\mathrm{mult}_{D_i}^{-1}(n_i).
$$
$\blacktriangleright$ Supposons que $\vert I\vert\geq 1$. S'il existe $i\in I$ tel que $n_{i}\not\in\{0,1\}$, il d\'ecoule de la propri\'et\'e \ref{mult1}/4 que~:
$$  
\pi_{|n|,Y}(A_n)=\emptyset.
$$
De m\^eme, comme $Y\rightarrow R$ est suppos\'ee lisse, les $\mathsf{D}_{i}$, pour $i\in I$, sont deux \`a deux disjoints. En particulier, pour tout $1\leq i\not=j\leq q$,
$$
\mathrm{mult}_{D_{i}}^{-1}(1)\cap\mathrm{mult}_{D_{i}}^{-1}(1)=\emptyset.
$$
Enfin, comme $\mathrm{mult}_{\mathsf{Y}}(y)=1$ pour tout $y\in\G(Y)$, il est clair que~:
$$
\cap_{i\in I}^q\mathrm{mult}_{D_{i}}^{-1}(0)=\emptyset,
$$
si $\vert I\vert=r\geq 1$. En outre, $A_n$ est vide si et seulement si l'une de ces conditions est v\'erifi\'ee, puisque $Y$ est lisse sur $R$. Par suite, $A_n=\emptyset$ si et seulement si $n\not\in\mathfrak{P}_D$. 

\bigskip

$\blacktriangleright$~Supposons d\'esormais que $n\in\mathfrak{P}_D$. Comme par hypoth\`ese $Y\rightarrow R$ est lisse, le $k$-morphisme de sch\'emas canonique $\theta_{0^n}~: \G_{\vert n\vert}(\!Y)\rightarrow\G_0(\!Y)$ induit une application surjective~:
$$
\pi_{\vert n\vert,\!Y}(A_{n})\rightarrow \pi_{0,\!Y}(A_{n}).
$$
Le th\'eor\`eme 4.2.3 de \cite{s} assure qu'il nous suffit de prouver une propri\'et\'e analogue sur les fibres de cette application pour d\'emontrer le r\'esultat. La question \'etant locale, on peut ramener l'\'etude au cas d'un $R$-sch\'ema affine $Y\rightarrow R$ lisse purement de dimension relative $d$, muni d'un $R$-morphisme de sch\'emas \'etale $Y\rightarrow \mathbf{A}_R^d$ (ou, ce qui revient au m\^eme, d'un syst\`eme de $d$ coordonn\'ees locales $x_1,\ldots,x_d$ tel que $D$ soit d\'efini par l'annulation d'un mon\^ome en les $x_i$. Soit $d_0$ un entier naturel non nul, inf\'erieur ou égal \`a $d$, et tel que l'\'equation de $D$ soit $x_{\ell_1}\ldots x_{\ell_{d_0}}$. On peut supposer que $J\subset\{1,\ldots,{d_0}\}$. 

Si, $x\in A_n$, la fibre au-dessus de $\pi_{0,Y}(x)$ par la restriction du morphisme de transition $\pi_{\vert n\vert,Y}(A_n)\rightarrow \pi_{0,Y}(A_n)$ est isomorphe \`a $\mathbf{G}_m^{\vert J\vert}\times \mathbf{L}^{\vert n\vert d-\vert n\vert}$. En particulier~:
$$
[\pi_{\vert n\vert, Y}(A_n)]\mathbf{L}^{-\vert n\vert d}=[\mathsf{D}_J^\circ]\times(\mathbf{L}-\mathbf{1})^{\vert J\vert}\times \mathbf{L}^{-\vert n\vert}.
$$ 
\end{proof}

\begin{prop}
\label{rationalite zeta dcn}
Soit $Y$ une $R$-vari\'et\'e lisse, purement de dimension relative $d$. Soit $D\hookrightarrow Y$ un diviseur réduit de $Y$. Si $D$ et $tD$ sont \`a croisements normaux stricts, alors la fonction Z\^eta $Z_D(T)$ appartient au sous-anneau de $\!M_k[[T]]$ engendr\'e par $\!M_k[T]$ et les \'el\'ements $(\mathbf{1}-\mathbf{L}^{-1}T^{i})^{-1}$ avec $i\in\mathbf{N}^m$. Plus pr\'ecis\'ement~:

\bigskip
			
$\vartriangleright$ si $\vert I\vert=0$~:
$$
Z_D(T)=\sum_{J\subset\{1,\ldots,m\}}[\mathsf{D}_J^\circ].\Bigl(\prod_{j\in J}\frac{(\mathbf{L}-\mathbf{1})\mathbf{L}^{-1}T_j}{\mathbf{1}-\mathbf{L}^{-1}T_j}\Bigr).
$$
$\vartriangleright$ Si $\vert I\vert\geq 1$ et $D\not\subset \mathsf{Y}$ et $\mathsf{Y}\not\subset D$~:
$$
Z_D(T)=\frac{\mathbf{L}-\mathbf{1}}{\mathbf{L}}\sum_{\ell\in I}\sum_{J\subset\{1,\ldots,m\}\backslash I}[\mathsf{D}_{J\cup\{\ell\}}^\circ]\Bigl(\prod_{j\in J}\frac{(\mathbf{L}-\mathbf{1})\mathbf{L}^{-1}T_j}{\mathbf{1}-\mathbf{L}^{-1}T_j}\Bigr)T_\ell
$$
$$
+\sum_{J\subset\{1,\ldots,m\}\backslash I}[\mathsf{D}_J^\circ].\Bigl(\prod_{j\in J}\frac{(\mathbf{L}-\mathbf{1})\mathbf{L}^{-1}T_j}{\mathbf{1}-\mathbf{L}^{-1}T_j}\Bigr).
$$
$\vartriangleright$ Si $\vert I\vert=r$ et $m>\vert I\vert$~:
$$
Z_D(T)=\frac{\mathbf{L}-\mathbf{1}}{\mathbf{L}}\sum_{\ell\in I}\sum_{J\subset\{1,\ldots,m\}\backslash I}[\mathsf{D}_{J\cup\{\ell\}}].\Bigl(\prod_{j\in J}\frac{(\mathbf{L}-\mathbf{1})\mathbf{L}^{-1}T_j}{\mathbf{1}-\mathbf{L}^{-1}T_j}\Bigr)T_\ell.
$$
$\vartriangleright$ Si $D=\mathsf{Y}$~:
$$
Z_D(T)=\frac{\mathbf{L}-\mathbf{1}}{\mathbf{L}}\sum_{i=1}^r\left[\mathsf{Y}_i\right]T_i.
$$
\end{prop}

\begin{proof}
Quitte \`a r\'eordonner les composantes irr\'eductibles de $D$, on peut supposer qu'il existe un entier $0\leq q\leq m$ tel que les composantes $D_i\hookrightarrow Y$ sont contenues dans $\mathsf{Y}$ si et seulement si $m-q+1\leq i\leq m$.

\bigskip

$\blacktriangleright$ Si $D=\mathsf{Y}$, la formule d\'ecoule directement du lemme \ref{calcul dcn}.

\bigskip

$\blacktriangleright$ Si $q=0$, il d\'ecoule de la d\'efinition et du lemme \ref{calcul dcn}~:
$$
Z_D(T)=\sum_{J\subset\{1,\ldots,m\}}\sum_{n\in(\mathbf{N}\backslash\{0\})^{\vert J\vert}}\left(\cap_{i\in J}\mathrm{mult}^{-1}_{D_i}(n_i)\right)\cap\left(\cap_{j\not\in J}\mathrm{mult}_{D_j}^{-1}(0)\right).
$$
Toujours gr\^ace au lemme \ref{calcul dcn}, on a plus pr\'ecis\'ement~:
$$
Z_D(T)=\sum_{J\subset\{1,\ldots,m\}}\left[D_J^\circ\right](\mathbf{L}-\mathbf{1})^{\vert J\vert}\sum_{n\in(\mathbf{N}\backslash\{0\})^{\vert J\vert}}\mathbf{L}^{-\vert n\vert}T^n.
$$
Le lemme \ref{calcul serie} implique alors que~:
$$
Z_D(T)=\sum_{J\subset\{1,\ldots,m\}}\left[D_J^\circ\right]\prod_{j\in J}\frac{(\mathbf{L}-\mathbf{1})\mathbf{L}^{-1}T_j}{\mathbf{1}-\mathbf{L}^{-1}T_j}.
$$

\bigskip

$\blacktriangleright$ Supposons que $D\not=\mathsf{Y}$. Le lemme \ref{calcul dcn} et l'additivit\'e de $[~]$ sur les constructibles assurent que~:
$$
Z_D(T)=\sum_{\ell=m-q+1}^m\sum_{J\subset\{1\ldots m-q\}}
$$
$$
\sum_{n\in\mathbf{N}^{m-q}}\left(\bigcap_{i\in J}\mathrm{mult}_{D_i}^{-1}(n_i)\right)\cap\left(\bigcap_{j\not\in J,j\in\{m-q+1,\ldots,m\}\backslash\{\ell\}}\mathrm{mult}_{D_j}^{-1}(0)\right)\cap \mathrm{mult}_{D_\ell}^{-1}(1)\mathbf{L}^{-(\vert n\vert+1)d}T^nT_\ell
$$
$$
 + 
$$
$$
(1-\delta_{q,r})\sum_{J\subset\{1,\ldots,m-q\}}\sum_{n\in\mathbf{N}^{m-q}}\left(\bigcap_{i\in J}\mathrm{mult}_{D_i}^{-1}(n_i)\right)\cap\left(\bigcap_{j\not\in J,j\in\{m-q+1,\ldots,m\}}\mathrm{mult}_{D_j}^{-1}(0)\right)\mathbf{L}^{-\vert n\vert}T^n, 
$$
o\`u $\delta_{q,r}$ vaut 1 si $q=r$, et 0 sinon. Les formules d\'ecoulent alors du lemme \ref{calcul dcn} par les m\^emes calculs que ceux du cas $q=0$.
\end{proof}

\begin{lem}
\label{calcul dcn x}
Soient $n\in\mathbf{N}^m$ et $x\in\mathsf{D}(k)$, alors on a dans $\!M_k$ la relation~:
$$
\left[\pi_{|n|,Y}\left((\cap_{i=1}^m\mathrm{mult}_{D_i}^{-1}(n_i))\cap \pi_{0,Y}^{-1}(x)\right)\right]=(\mathbf{L}-\mathbf{1})^{\vert I_x\vert}\times\mathbf{L}^{\vert n\vert(d-1)},
$$
si $n\in\mathfrak{P}_D$ et si $J= I_x$. Sinon,
$$
\left[\pi_{|n|,Y}\left((\cap_{i=1}^m\mathrm{mult}_{D_i}^{-1}(n_i))\cap \pi_{0,Y}^{-1}(x)\right)\right]=0.
$$
\end{lem}

\begin{proof}
Posons $A_{n,x}:=(\cap_{i=1}^m\mathrm{mult}_{D_i}^{-1}(n_i))\cap \pi_{0,Y}^{-1}(x)$. Comme par hypoth\`ese $Y\rightarrow R$ est lisse, le $k$-morphisme de sch\'emas canonique $\theta_{0^n}~: \G_{\vert n\vert}(\!Y)\rightarrow\G_0(\!Y)$ induit une application surjective~:
$$
\pi_{\vert n\vert,\!Y}(A_{n,x})\rightarrow \pi_{0,\!Y}(A_{n,x})=\{x\},
$$
dont les fibres, si $A_{n,x}\not=\emptyset$, sont isomorphes \`a $\mathbf{G}_{m,\kappa(x)}^{\vert I_x\vert}\times_k\mathbf{A}_{\kappa(x)}^{\vert n\vert(d-1)}$. Le lemme d\'ecoule donc du th\'eor\`eme 4.2.3 de \cite{s}.

\end{proof}

\begin{prop}
\label{rationalite zeta dcn x}
Soit $Y$ une $R$-vari\'et\'e lisse, purement de dimension relative $d$. Soient $D\hookrightarrow Y$ un diviseur réduit de $Y$ et $x\in\mathsf{D}(k)$. Si $D$ et $tD$ sont \`a croisements normaux stricts, alors la fonction Z\^eta $Z_{D,x}(T)$ appartient au sous-anneau de $\!M_k[[T]]$ engendr\'e par $\!M_k[T]$ et les \'el\'ements $(\mathbf{1}-\mathbf{L}^{-1}T^{i})^{-1}$ avec $i\in\mathbf{N}^m$. Plus pr\'ecis\'ement,

\smallskip

$\rhd$~Si $D\not=\mathsf{Y}$ et si $x$ n'appartient \`a aucune composante de $\mathsf{Y}$ contenue dans $\mathsf{D}$~: 
$$
Z_{D,x}(T)=\prod_{i\in I_x}\frac{(\mathbf{L}-\mathbf{1})\mathbf{L}^{-1}T_i}{\mathbf{1}-\mathbf{L}^{-1}T_i}.
$$
$\rhd$~Si $D\not=\mathsf{Y}$ et si $x$ appartient \`a une composante $\mathsf{Y}_\ell$ de $\mathsf{Y}$ contenue dans $\mathsf{D}$~: 
$$
Z_{D,x}(T)=\frac{\mathbf{L}-\mathbf{1}}{\mathbf{L}}T_\ell\prod_{i\in I_x,i\not=\ell}\frac{(\mathbf{L}-\mathbf{1})\mathbf{L}^{-1}T_i}{\mathbf{1}-\mathbf{L}^{-1}T_i}.
$$
$\rhd$~Si $D=\mathsf{Y}$~:
$$
Z_{D,x}(T)=\frac{\mathbf{L}-\mathbf{1}}{\mathbf{L}}T_\ell.
$$

\end{prop}

\begin{proof}
Soit $x\in\mathsf{D}(k)$. Posons $I:=I_x$. Par d\'efinition,
$$
Z_{D,x}(T)=\sum_{n\in\mathbf{N}^m}\left[\pi_{|n|,Y}\left((\cap_{i=1}^m\mathrm{mult}_{D_i}^{-1}(n_i))\cap \pi_{0,Y}^{-1}(x)\right)\right]\mathbf{L}^{-\vert n\vert}T^n
$$
Il d\'ecoule du lemme \ref{calcul dcn x} que~:
$$
Z_{D,x}(T)=\sum_{n\in(\mathbf{N}\backslash\{0\})^{\vert I\vert}}\left[\pi_{|n|,Y}\left(\left(\cap_{i\in I}\mathrm{mult}_{D_i}^{-1}(n_i)\right)\cap\left(\cap_{i\not\in I}\mathrm{mult}_{D_i}^{-1}(0)\right)\cap \pi_{0,Y}^{-1}(x)\right)\right]\mathbf{L}^{-\vert n\vert}T^n.
$$
Si $x$ n'appartient \`a aucune composante de $D$ contenue dans $\mathsf{Y}$, toujours gr\^ace au lemme \ref{calcul dcn x},
$$
Z_{D,x}(T)=\sum_{n\in(\mathbf{N}\backslash\{0\})^{\vert I\vert}}(\mathbf{L}-\mathbf{1})^{\vert I\vert}\mathbf{L}^{-\vert n\vert}T^n.
$$
Sinon, il existe une unique composante de $\mathsf{Y}$ contenue dans $D$ et contenant $x$ et donc~:
$$
Z_{D,x}(T)=\frac{\mathbf{L}-\mathbf{1}}{\mathbf{L}}T_\ell\sum_{n\in(\mathbf{N}\backslash\{0\})^{\vert I\vert-1}}(\mathbf{L}-\mathbf{1})^{\vert I\vert-1}\mathbf{L}^{-\vert n\vert}T^n.
$$
Les formules d\'ecoulent alors directement du lemme \ref{calcul serie}.

\end{proof}

\begin{prop}
\label{P calcul}
Soit $X$ une $R$-vari\'et\'e lisse, purement de dimension relative $d$. Soit $E\hookrightarrow X$ une sous-$R$-vari\'et\'e ferm\'ee de $X$, de codimension au moins 1, telle que $tE$ contient le lieu singulier de $X$. Alors l'int\'egrale~:
\begin{eqnarray}
\int_{\G(X)}\mathbf{L}^{-\mathrm{mult}_{E}}d\mu_{X}\in\widehat{\!M_k}.
\end{eqnarray}
appartient au sous-anneau engendr\'e par $\overline{\!M_k}$ et les \'el\'ements $(\mathbf{L}^a-\mathbf{1})^{-1}$ pour $a\in\mathbf{N}\backslash\{0\}$. 
\end{prop}

\begin{proof}
La fonction $x\in\G(X)\mapsto\mathrm{mult}_E(x)$ est exponentiellement int\'egrable car $E$ est de codimension au moins 1 dans $X$. Soit $h~:Y\rightarrow X$ une r\'esolution de N\'eron plong\'ee de $(X,tE)$. Le th\'eor\`eme de changement de variables \ref{tcvg} assure alors que~:
$$
I_E:=\int_{\G(X)}\mathbf{L}^{-\mathrm{mult}_E}d\mu_X=\int_{\G(Y)}\mathbf{L}^{-\mathrm{mult}_{h^{-1}(E)}-\ord_t(\mathrm{Jac})_h}d\mu_Y.
$$
Comme $X$ est lisse, $\ord_t(\mathrm{Jac})^{-1}(\infty)$ est un diviseur de $Y$ contenu dans $D:=h^{-1}(E)$ qui est un diviseur de $Y$ \`a croisements normaux stricts, par d\'efinition des r\'esolutions de N\'eron plong\'ees. Soient $D_i$, pour $1\leq i\leq m$, les composantes irr\'eductibles de $D$. Par additivit\'e de la mesure $\mu_Y$ (cf. lemme 7.0.7 de \cite{s}), on peut supposer que $Y$ est affine, que les $a_i$ sont les multiplicit\'es de $D_i$ dans $D$ et $b_i$ celles de $D_i$ dans $\ord_t(\mathrm{Jac})^{-1}(\infty)$. On a donc~:
$$
I_E=\sum_{n\in\mathbf{N}}\mu_Y\left(\mathrm{mult}_{h^{-1}(E)}+\ord_t(\mathrm{Jac})_h\right)^{-1}(n)\mathbf{L}^{-(n+1)d}d\mu_Y.
$$
Le lemme \ref{mult3} assure que~:
$$
I_E=\mathbf{L}^{-d}\sum_{\ell\in\mathbf{N}^m}\left[\pi_{\vert \ell\vert,Y}\left(\cap_{i=1}^m\mathrm{mult}_{D_i}(\ell_i)\right)\right]\mathbf{L}^{-(\sum_{i=1}^{m}(a_i+b_i)\ell_i)d}.
$$
Or, par d\'efinition,
$$
Z_{D_{\mathrm{red}}}(T):=\sum_{n\in\mathbf{N}^m}\left[\pi_{\vert n\vert,Y}\left(\cap_{i=1}^m\mathrm{mult}_{D_i}(n_i)\right)\right]\mathbf{L}^{-\vert n\vert d}T^n.
$$
Si l'on d\'esigne par $\iota~:\!M_k\rightarrow\overline{\!M_k}$ le morphisme d'anneaux canonique, on a la relation~:
$$
\iota\left(\mathbf{L}^{-d}Z_{D_{\mathrm{red}}}\left(\mathbf{L}^{-(a_1+b_1-1)d},\mathbf{L}^{-(a_2+b_2-1)d},\ldots,\mathbf{L}^{-(a_m+b_m-1)d}\right)\right)=I_E.
$$
La proposition d\'ecoule alors de l'assertion de la proposition \ref{rationalite zeta dcn}.
\end{proof}

\subsection{La fonction Z\^eta motivique associ\'ee \`a une famille de fonctions d'une vari\'et\'e lisse}

Soit $X$ une $R$-vari\'et\'e lisse, purement de dimension relative $d$. Soit $f:=(f_1,\ldots,f_\ell)~:X\rightarrow \mathbf{A}^\ell_R$ un $R$-morphisme de sch\'emas. Notons $E_i:=V(f_i)\hookrightarrow X$ pour tout $1\leq i\leq \ell$. On d\'efinit alors la fonction Z\^eta associ\'ee \`a $f$ comme la s\'erie de ${\!M_k}[[T_1,\ldots,T_\ell]]$~:
\begin{eqnarray}
Z_f(T):=\sum_{n\in\mathbf{N}^\ell}\left[\pi_{\vert n\vert,X}\big(\cap_{i=1}^\ell\mathrm{mult}_{E_i}^{-1}(n_i)\big)\right]\mathbf{L}^{-\vert n\vert d}T^n,
\end{eqnarray}
avec les notations pr\'ec\'edentes. Soit $E$ le diviseur effectif principal de $Y$ d\'efini par le produit des $f_i$, pour $1\leq i\leq \ell$. Si $x\in\mathsf{E}(k)$, la fonction Z\^eta de $f$ de \textit{support $x$} est la s\'erie~:
$$
Z_{f,x}(T):=\sum_{n\in\mathbf{N}^\ell}\left[\pi_{\vert n\vert,X}\left(\big(\cap_{i=1}^\ell\mathrm{mult}_{E_i}^{-1}(n_i)\big)\cap\pi_{0,X}^{-1}(x)\right)\right]\mathbf{L}^{-\vert n\vert d}T^n.
$$

\begin{prop}
\label{Premier calcul}
Soit $Y$ une $R$-vari\'et\'e lisse, purement de dimension relative $d$. Soit $f:=(f_1,\ldots,f_\ell)~:Y\rightarrow \mathbf{A}^\ell_R$ un $R$-morphisme de sch\'emas. Alors les fonctions Z\^eta $Z_f(T)$ et $Z_{f,x}(T)$ (avec $x\in \mathsf{E}(k)$) dans $\!M_k[[T_1,\ldots,T_\ell]]$ appartiennent au sous-anneau engendr\'e par $\!M_k[T_1,\ldots,T_\ell]$ et les \'el\'ements $(\mathbf{1}-\mathbf{L}^{-b}T^{a})^{-1}$ avec $a\in\mathbf{N}^\ell$ et $b\in\mathbf{N}$.
\end{prop}

\begin{proof}
Les preuves \'etant identiques, nous nous contentons du cas de $Z_f(T)$. Soit $h~: Y\rightarrow X$ une r\'esolution de N\'eron plong\'ee (cf. \ref{desingularisation Neron}) de $(X,tE)$. La lissit\'e de $X$ et le lemme \ref{dcn} assurent que $\ord(\mathrm{Jac})_h^{-1}(\infty)$ est un diviseur \`a croisements normaux stricts contenu dans $C:=h^{-1}(E)$. Soient $C_i$, pour $1\leq i\leq m$, les composantes irr\'eductibles de $C$. On note $\nu_i-1$ la multiplicit\'e de $C_i$ dans $\ord(\mathrm{Jac})_h^{-1}(\infty)$. De la m\^eme mani\`ere, $D_j:=h^{-1}(E_j)$, pour $1\leq j\leq \ell$, est un diviseur \`a croisements normaux stricts de $Y$ et l'on d\'esigne par $a_{ij}$ la multiplicit\'e de $C_i$ dans $D_j$, pour tout $i$ et tout $j$. Commen\c{c}ons par remarquer que, pour tout $1\leq i\leq \ell$ et tout $n\in\mathbf{N}^\ell$, comme $\ord(\mathrm{Jac})_h^{-1}(\infty)\subset D$,
$$
\cap_{i=1}^\ell\mathrm{mult}_{D_i}^{-1}(n_i)=\bigsqcup_{e\in\mathbf{N}}\left(\left(\cap_{i=1}^\ell\mathrm{mult}_{D_i}^{-1}(n_i)\right)\cap\ord_t(\mathrm{Jac})_h^{-1}(e)\right),
$$
et que cette r\'eunion, par quasi-compacit\'e de la topologie constructible, peut \^etre suppos\'ee finie. Soit
$$
A_{n,e}:=\left(\cap_{i=1}^\ell\mathrm{mult}_{D_i}^{-1}(n_i)\right)\cap\ord_t(\mathrm{Jac})_h^{-1}(e).
$$
Il d\'ecoule du lemme 8.1.3 de \cite{s}, qui est une sorte de formule de changement de variables dans $\!M_k$, de l'additivit\'e de $\left[~\right]$ sur les constructibles et du lemme \ref{mult2} que~:
$$
\left[\cap_{i=1}^\ell\mathrm{mult}_{E_i}^{-1}(n_i)\right]=\sum_{e\in\mathbf{N}}\left[ A_{n,e}\right]\mathbf{L}^{-e}.
$$
Soit alors $Z(T,U)$ la s\'erie formelle de $\!M_k[[T,U]]$ d\'efinie par~:
$$
Z(T,U)=\sum_{(n,e)\in\mathbf{N}^\ell\times\mathbf{N}}\left[A_{n,e}\right]T^nU^e.
$$
Il d\'ecoule de cette d\'efinition que~:
$$
Z_f(T)=Z(T,\mathbf{L}^{-1}).
$$
Le lemme \ref{mult3} assure que~:
$$
Z(T,U)=\sum_{s\in\mathbf{N}^m}\left[\cap_{i=1}^m\mathrm{mult}^{-1}_{C_i}(s_i)\right]U^{\sum_{i=1}^m(\nu_i-1)}T_1^{\sum_{i=1}^ma_{1i}s_i}\ldots T_\ell^{\sum_{i=1}^ma_{\ell i}s_i}.
$$
La relation~:
$$
Z(T,U)=Z_{C_{\mathrm{red}}}(U^{\nu_1-1}T_1^{a_{11}s_1}T_2^{a_{21}s_1}\ldots T_\ell^{a_{\ell 1}s_1},\ldots,U^{\nu_m-1}T_1^{a_{1m}s_m}T_2^{a_{2m}s_m}\ldots T_\ell^{a_{\ell m}s_m})
$$
implique le r\'esultat gr\^ace \`a la proposition \ref{rationalite zeta dcn}.
\end{proof}

\begin{rem}
L'\'egalit\'e 
$$
Z_f(T)=Z_{C_{\mathrm{red}}}(\mathbf{L}^{1-\nu_1}\prod T_j^{a_{j1}s_1},\ldots,\mathbf{L}^{1-\nu_1}\prod T_j^{a_{jm}s_m})
$$ 
donne des formules explicites pour les fonctions Z\^eta, gr\^ace \`a la proposition \ref{rationalite zeta dcn}.
\end{rem}

\subsection{La fonction Z\^eta associ\'ee \`a une famille de fonctions d'une vari\'et\'e}

Soit $X$ une $R$-vari\'et\'e plate, purement de dimension relative $d$. Soit $f:=(f_1,\ldots,f_\ell)~:X\rightarrow \mathbf{A}^\ell_R$ un $R$-morphisme de sch\'emas. Notons $E_i:=V(f_i)\hookrightarrow X$ pour tout $1\leq i\leq \ell$. On d\'efinit alors la fonction Z\^eta associ\'ee \`a $f$ comme la s\'erie de ${\!M_k}[[T_1,\ldots,T_\ell]]$~:
\begin{eqnarray}
Z_f(T):=\sum_{n\in\mathbf{N}^\ell}\mu_X\left(\cap_{i=1}^\ell\mathrm{mult}_{E_i}^{-1}(n_i)\right)T^n,
\end{eqnarray}
avec les notations pr\'ec\'edentes. Soit $E$ le diviseur effectif principal de $X$ d\'efini par le produit des $f_i$, pour $1\leq i\leq \ell$. Si $x\in\mathsf{E}(k)$, la fonction Z\^eta de $f$ de \textit{support $x$} est la s\'erie~:
$$
Z_{f,x}(T):=\sum_{n\in\mathbf{N}^\ell}\mu_X\left(\big(\cap_{i=1}^\ell\mathrm{mult}_{E_i}^{-1}(n_i)\big)\cap\pi_{0,X}^{-1}(x)\right)T^n.
$$

\begin{prop}
\label{X singulier}
Soit $X$ une $R$-vari\'et\'e plate, purement de dimension $d+1$. Soit $f:=(f_1,\ldots,f_\ell)~:X\rightarrow \mathbf{A}^\ell_R$ un $R$-morphisme de sch\'emas. Soit $E$ le diviseur effectif principal de $Y$ d\'efini par le produit des $f_i$, pour $1\leq i\leq \ell$. Alors les fonctions Z\^eta $Z_f(T)$ et $Z_{f,x}(T)$ (avec $x\in \mathsf{E}(k)$) de $\widehat{\!M_k}[[T_1,\ldots,T_\ell]]$ appartiennent au sous-anneau engendr\'e par $\overline{\!M_k}[T_1,\ldots,T_\ell]$, les \'el\'ements $(\mathbf{1}-\mathbf{L}^{-b}T^{a})^{-1}$ et $(\mathbf{L}^i-\mathbf{1})^{-1}$, avec $a\not=0\in\mathbf{N}^\ell$, $b\in\mathbf{N}$ et $i\in\mathbf{N}\backslash\{0\}$.

\end{prop}

\begin{proof}
L'exemple \ref{ordjac} assure que les familles $(A_n)_{n\in\mathbf{N}^{l}}$ et $(B_n)_{n\in\mathbf{N}^{\ell}}$, d\'efinies par $A_n:=\cap_{i=1}^\ell\mathrm{mult}_{E_i}^{-1}(n_i)$ et $B_n:=A_n\cap \pi_{0,X}^{-1}(x)$, sont semi-alg\'ebriques. On conclut gr\^ace aux th\'eor\`emes \ref{thm3} et \ref{thm2}, en remarquant auparavant que~:
$$
\mu_X(A_n)=\int_{A_n}\mathbf{L}^{-0}d\mu_X
$$
et
$$
\mu_X(B_n)=\int_{B_n}\mathbf{L}^{-0}d\mu_X.
$$

\end{proof}

\begin{rem}
\label{looijenga}
Comme pour l'anneau $\!M_k$, si $S$ est un sch\'ema, on peut d\'efinir l'anneau $\!M_S$ \`a partir de la cat\'egorie des $S$-vari\'et\'es. Si $n\in\left(\mathbf{N}\backslash\{0\}\right)^\ell$, on peut d\'efinir la classe de $\cap_{i=1}^\ell\mathrm{mult}_{E_i}^{-1}(n_i)$ dans l'anneau $\!M_{(\mathbf{G}_{m,k}\times_k \mathsf{E})^\ell}:=\!M_{\mathbf{G}_{m,\mathsf{E}}^\ell}$. En particulier, on pourrait \'etudier des fonctions Z\^eta motiviques `` non triviales '' de $\!M_{\mathbf{G}_{m,\mathsf{E}}^\ell}$. C'est, notamment, ce que fait Looijenga dans \cite{lo}. La proposition 4.2 de \textit{loc. cit.} est un analogue de notre lemme \ref{calcul dcn} et le morphisme canonique (d'oubli)~:
$$
\!M_{\mathbf{G}_{m,\mathsf{E}}^\ell} \rightarrow \!M_k
$$
sp\'ecialise les fonctions Z\^eta de \textit{loc. cit.} en nos fonctions Z\^eta (cf. corollaire 4.3 de \textit{loc. cit.}). Il est important de noter que l'introduction des r\'esolutions de N\'eron plong\'ees donnent un sens aux hypoth\`eses de la proposition 4.2. 

Notons enfin, que dans le cas d'une vari\'et\'e singuli\`ere, l'analogue dans $\!M_{\mathbf{G}_{m,\mathsf{E}}^\ell}$ de la proposition \ref{X singulier}, demanderait de g\'en\'eraliser, \`a l'anneau $\!M_{\mathbf{G}_{m,\mathsf{E}}^\ell}$, les th\'eor\`emes \ref{thm3} et \ref{thm2}. Il ne nous semble pas qu'un tel \'enonc\'e soit d\'emontr\'e dans \cite{lo}.
\end{rem}

\bibliographystyle{amsplain}

\end{document}